\documentclass{svmult}
\usepackage{amssymb,amsmath,amsfonts}

\usepackage{mathtools}
\usepackage{bm}

\def\dd{\mathrm{d}}

\usepackage{tikz}
\usetikzlibrary{arrows,trees,decorations,external,positioning,backgrounds}

\usepackage{listings}
\usepackage{color}
\definecolor{dkgreen}{rgb}{0,0.6,0}
\definecolor{gray}{rgb}{0.5,0.5,0.5}
\definecolor{mauve}{rgb}{0.58,0,0.82}

\usepackage{float}

\graphicspath{{imgs/}} %Setting the graphicspath

\graphicspath{{imgs/}}
\makeatletter
\def\input@path{{imgs/}}
\makeatother

\begin{document}

\title{Towards Geometric Time Minimal Control without Legendre Condition 
and with Multiple Singular Extremals for Chemical Networks}

\titlerunning{Geometric Techniques to Optimize Chemical Reactions}

\author{Bernard Bonnard\inst{1}\and
J\'er\'emy Rouot\inst{2}}

\institute{
  Inria Sophia Antipolis and Institut de Math\'ematiques de Bourgogne, 9 avenue Savary, Dijon, 21078 France,\\
  \texttt{bernard.bonnard@u-bourgogne.fr}
\and
\texttt{jeremy.rouot@grenoble-inp.org}}

\maketitle

\abstract{
  This article deals with the problem of maximizing the production of a species 
  for a chemical network by controlling the temperature. 
  Under the so--called mass kinetics assumption the
  system can be modeled as a single--input control system using the
  Feinberg--Horn--Jackson graph associated to the reactions network. 
  Thanks to Pontryagin's Maximum Principle, the candidates as minimizers can be found
  among extremal curves, solutions of a (non smooth) Hamiltonian dynamics 
  and the problem can be stated as a time minimal control problem with a terminal target
  of codimension one. Using geometric control and singularity theory the time minimal
  syntheses (closed loop optimal control) can be classified near the terminal manifold 
  under generic conditions.
  In this article, we focus to the case where the generalized Legendre-Clebsch condition
  is not satisfied, which paves the road to complicated syntheses with several singular
  arcs.
  In particular, it is related to the situation for a weakly reversible network like
  the McKeithan scheme of two reactions:
$
\tikz[baseline={(a.base)},node distance=6mm] {
  \node (a) {T+M};
    \node[right=of a] (b) {A};
    \node[right=of b] (c) {B};
\draw (a.base east) edge[left,->] node[above] {} (b.base west) (b.base east) edge[->] node[above] {} (c.base west);
    \draw (a) edge[bend right=40,<-] node[below right] {} (b);
   \draw (a.south) edge[bend right=60,<-] node[below right] {} (c);
  }
  $.
}



%% pour eviter le numero du chapitre devant dans numerotations
\renewcommand{\thesection}{\arabic{section}}
%% profondeur des numerotaitons
\setcounter{secnumdepth}{5}

\section{Introduction}
\label{sec:intro}
The optimization of the production is an important problem in chemical and biological engineering
and the control can be either the temperature (batch or closed reactor) or by feeding the
reactor (semi-batch or open case).
In this article we shall concentrate on the first case. Moreover we assume that the
dynamics is modeled using the mass action kinetics assumptions and hence given
at constant temperature $T$ by a polynomial system based only on the Feinberg--Horn--Jackson
graph associated to the chemical network.
Also in the 70's those researchers obtained (under the so--called zero deficiency assumption)
in a series of seminal articles \cite{feinberg1972,horn1972} a complete description of the 
dynamics, at constant temperature.
If this dynamics is well understood, in the optimal problem the temperature is not constant
and the analysis becomes very intricate.
Thanks to the Pontryagin Maximum Principle \cite{pontryagin1962} candidates as minimizers
can be found among extremals solutions of a ({\it non smooth}) Hamiltonian dynamics and 
optimal solutions are concatenation of bang arcs, with minimum and maximum temperature, and 
the so--called singular arcs, defined as a solution of a {\it smooth} Hamiltonian 
{\it constrained} dynamics \cite{bonnard2003}. Moreover maximizing the production of one species
during the batch can be restated as producing a fixed amount of this species while minimizing
the batch duration.
In this frame, the problem is a time minimal control problem, 
with a terminal target manifold
of {\it codimension one}. 

For applications, the time optimal control has to be computed as a closed loop feedback and this
leads to the problem of computing the time minimal {\it synthesis}, for a single--input
control system. At the end of the 80's, geometric optimal control has produced an important 
literature to compute the time minimal syntheses for the {\it fixed end point case}, 
where the initial
point is localized near the terminal point. 
This was done using the {\it Lie algebraic structure} of the control system, mainly for single--input
(smooth) system, where the control appears linearly, some seminal references are
\cite{sussmann1990,schattler1988,boscain2004} either in general context or in view of application formed by a
sequence of two irreversible reactions: $A\rightarrow B \rightarrow C$ and in relation with an
industrial project \cite{bonnard1998}.

Our aim is to extend this work to more complicated reaction schemes and to deal in particular
with weakly reversible chemical schemes. More precisely we shall concentrate on a McKeithan
type scheme of the form 
$
\tikz[baseline={(a.base)},node distance=6mm] {
  \node (a) {T+M};
    \node[right=of a] (b) {A};
    \node[right=of b] (c) {B};
\draw (a.base east) edge[left,->] node[above] {} (b.base west) (b.base east) edge[->] node[above] {} (c.base west);
    \draw (a) edge[bend right=40,<-] node[below right] {} (b);
   \draw (a.south) edge[bend right=60,<-] node[below right] {} (c);
  }
  $
  assuming that the coefficient governing the dynamics are given by Arrhenius law.
  This scheme was already studied in the context of control theory using stabilization techniques
  with ``{\it feeding}'' types control \cite{sontag2001} 
  and ad hoc {\it observer design} \cite{chavez2002a}.
  In our case, this network is a test bed case for our very general approach.

  The key point for this extension is the analysis of singular trajectories and their role
  in the synthesis. This is connected with an important question and the need to extend the
  standard synthesis related to the turnpike phenomenon 
  \cite{sussmann1987a} to deal with
  cases, where the 
  strict Legendre--Clebsch condition is not satisfied, 
  a situation encountered in a recent
  application in MRI \cite{cots2018}.

  The organization of this article is the following. In Section \ref{sec2}, we recall
  briefly the Feinberg--Horn--Jackson theory to model the dynamics of chemical networks
  of constant temperature and the properties of the dynamics under the zero deficiency
  assumption \cite{feinberg1972,horn1972}, which can be applied to the McKeithan scheme.
  The stability properties are recalled and can be applied to control stabilization
  and observer design \cite{sontag2001,chavez2002a}.
  In Section \ref{sec3}, we present the fundamental results of the time minimal control problem,
  which are relevant to our study: Pontryagin Maximum Principle \cite{pontryagin1962}, regular
  and singular extremals \cite{bonnard2003},\cite{kupka1987}. The general turnpike theorem
  \cite{sussmann1987a} is recalled and extensions are presented in relation with the 
  problem with terminal manifold of codimension one and when the 
  strict Legendre--Clebsch condition is not
  satisfied. 
  The concept of conjugate and focal points is introduced based on \cite{bb1993}.
  In Section \ref{sec4}, we analyze the time minimal control problem for a two
  reactions McKeithan scheme. 
  To compute the time minimal syntheses, we present techniques and results from 
  \cite{bonnard1998,bonnard1997}, which have to be extended to analyze the problem.
  The computational complexity of the problem is discussed and symbolic computations are
  presented to cope with this complexity.

  \section{Mathematical model and stability properties of the McKeithan network}
  \label{sec2}
  \subsection{Mass action kinetics networks and dynamics using the Feinberg--Horn--Jackson
  graph (\cite{feinberg1972},\cite{horn1972})}

  We consider a set of $m$ chemical species $\{X_1,\ldots,X_m\}$ and the state of the dynamics
  is the vector $c=(c_1,\ldots,c_m)^\intercal \in \mathbb{R}_{\ge 0}^m$ representing the 
  molar concentration.
  Let ${\cal R}$ be a set of reactions, each reaction being denoted by $y\rightarrow y'$
  and of the form 
  \[
    \sum_{i=1}^m \alpha_i X_i \longrightarrow \sum_{i=1}^m \beta_i X_i,
  \]
  where $\alpha_i, \beta_i$ are the {\it stoichiometric coefficients}
  	and the vectors $y=(\alpha_1,\ldots, \alpha_m)^\intercal$
  and $y'=(\beta_1,\ldots, \beta_m)^\intercal$ are the vertices of the so--called
  {\it Feinberg--Horn--Jackson oriented graph} associated to the network, edges being 
  oriented according to $y\rightarrow y'$. Each reaction is characterized by a
  {\it reaction rate} $K(y\rightarrow y')$ and the system is said {\it simple} (or
  mass kinetics) if the rate of the reaction is of the form:
  \begin{equation}
    \begin{aligned}
    &K(y\rightarrow y') = k(T)\, c^y,
\\
&c^y = c_1^{\alpha_1} \ldots c_m^{\alpha_m}
    \end{aligned}
    \label{eq:constant-reaction}
  \end{equation}
  and 
  \[
    k(T) = A\, e^{-E/(RT)}
  \]
  is the {\it Arrhenius law}, $A$ is the exponential factor, $E$ is the activation
  energy, both depending on the reaction, $R$ is the gas constant and $T$ is the 
  temperature. Note that different rate formulae can be used to deal in particular 
  with biomedical systems (see for instance \cite{sontag2001}).
  The dynamics of the system, taking into account the whole network is:
  \begin{equation}
    \dot c(t) = f(c(t),T) = \sum_{y\rightarrow y'} K(y\rightarrow y')\, (y'-y).
    \label{eq:dyn-gen}
  \end{equation}

  \subsection{More explicit representation of the dynamics}
  \begin{definition}
    The stoichiometric subspace is 
    $S\coloneqq \mathrm{span}\{y-y';\ y\rightarrow y' \in {\cal R}\}$
    and the sets $(c(0) + S) \cap \mathbb{R}_{\ge 0}^m$ are called the (strictly if $>0$)
    positive stoichiometric compatibility classes.
    \label{def:stochio-subspace}
  \end{definition}
  From \cite{anderson2008} we have.
  \begin{lemma}
    Let $c(t)$ be a solution of \eqref{eq:dyn-gen} with initial condition $c(0) \in
    \mathbb{R}_{\ge 0}^m$. Then $c(t)$ belongs for all $t\ge 0$ to the strictly 
    positive compatibility class $(c(0)+S) \cap \mathbb{R}_{\ge 0}^m$.
    \label{lem:invariant-set}
  \end{lemma}
  \begin{definition}
    Having labeled the set of vertices by $i=1,\ldots, n$, with corresponding 
    stoichiometric vector $(y_1,\ldots,y_n)$, the complex matrix is 
    $Y\coloneqq(y_1,\ldots, y_n)$.
    The incidence connectivity matrix $A\coloneqq (a_{ij})$
    contains the Arrhenius coefficients $k_i$ of the reactions using the rule:
    $k_1=a_{21}$ indicates a reaction with kinetics constant $k_1$ from the first
    node to the second, that is $y_1 \underset{k_1}{\rightarrow} y_2$.
    \label{def:complex-matrix}
  \end{definition}

  With the mass kinetics assumption, the dynamics can be expressed as
  \[
    \dot c(t) = f(c(t),T) = Y \tilde A c^Y,
  \]
  where $\tilde A$ is the Laplacian matrix in graph theory defined by
  \begin{equation}
    \tilde A  = A - \text{diag}
    \left(
      \sum_{i=1}^n a_{i1},\ldots, \sum_{i=1}^n a_{in}
    \right)
    \label{eq:laplacian-matrix}
  \end{equation}
  and we denote 
  \[
    c^Y = (c^{y_1},\ldots,c^{y_n})^\intercal.
  \]

  \subsection{The McKeithan scheme (\cite{keithan2001,sontag2001})}
  It is given by the reaction scheme:
\[
\tikz[baseline={(a.base)},node distance=6mm] {
\node (a) {$T+M$};
\node[right=of a] (b) {$C_0$};
\node[right=of b] (c) {$C_1$};
\node[right=of c] (d) {$\hdots$};
\node[right=of d] (e) {$C_N$};
  \draw (a) edge[->] node[above] {$k_1$} (b) (b) edge[->] node[above] {$k_{p,0}$} (c) (c) edge[->] node[above] {$k_{p,1}$} (d) (d) edge[->] node[above] {$k_{p,N}$} (e);
  \draw (a) edge[bend right=30,<-] node[below right] {$k_{-1,0}$} (b);
  \draw (a) edge[bend right=53,<-] node[below] {$k_{-1,1}$} (c);
  \draw (a.south) edge[bend right=64,<-] node[below] {$k_{-1,N}$} (e);
}.
\]
The matrix $Y$ is given by
  \begin{equation*}
  Y =
  \begin{pmatrix}
    1 & 0 & \hdotsfor{2}& 0
    \\
    1 & 0 & & & \cdot
    \\[-1.5ex]
    0 & 1 & \ddots & & \vdots
    \\[-1.5ex]
    \vdots & \ddots & \ddots & \ddots & \vdots
    \\[-1.5ex]
   \vdots & & \ddots & 1 & 0
    \\
    0 & \hdotsfor{2} & 0 & 1
  \end{pmatrix}
 \end{equation*}
and the matrix $A=(a_{ij})$ is defined by $a_{21}=k_1,\ a_{1i}=k_{-1,i-2},\ i=2,\ldots,m$
$(m=N+2)$, $a_{i,i-1}=k_{p,i-3},\ i=3,\ldots,m$ and all others $a_{ij}$ are zero.

The stoichiometric subspace is defined by: 
\[{\{ c: T+C_0+\ldots+C_N = M+C_0+\ldots+C_N=0\}}\]
and we note: ${\delta_1=T+C_0+\ldots +C_N} \text{ and } \delta_2=M+C_0+\ldots +C_N$ 
the constants associated to first integrals of the dynamics.

Consider the special case $N=2$ so that the reaction scheme is denoted by 
$
\tikz[baseline={(a.base)},node distance=6mm] {
  \node (a) {T+M};
    \node[right=of a] (b) {A};
    \node[right=of b] (c) {B};
\draw (a.base east) edge[left,->] node[above] {{$\scriptstyle k_1$}} (b.base west) (b.base east) edge[->] node[above] {{$\scriptstyle k_2$}} (c.base west);
    \draw (a) edge[bend right=40,<-] node[below right] {{$\scriptstyle k_3$}} (b);
   \draw (a.south) edge[bend right=60,<-] node[below right] {{$\scriptstyle k_4$}} (c);
  }
  $,
  and restricting to the stoichiometric class ($\delta_1$ and $\delta_2$ being fixed), one 
  gets with: $x\coloneqq [A], y\coloneqq [B]$, $[\cdot]$ denoting the respective concentrations, 
  that the
  dynamics is described by the equations
  \begin{equation}
    \begin{aligned}
      &\dot x = k_1(\delta_1 x-y)(\delta_2 -x -y)-(k_2+k_3)x
      \\
      &\dot y = k_2x-k_4 y.
    \end{aligned}
    \label{eq:dyn}
  \end{equation}

  \begin{definition}
    The deficiency of the network is: $\delta=n-l-s$, where $n$ is the number of vertices, 
    $l$ is the number of connected components and $s$ is the dimension of the 
    stoichiometric subspace. The network is called strongly connected if for each pair $(i,j)$
    of vertices such that there exists an oriented path joining $i$ to $j$ there exists 
    a path joining $j$ to $i$.
  \end{definition}

    Using \cite{feinberg1972}, refined by \cite{anderson2008,sontag2001}, one has
    the following result.
    \begin{theorem}
      The graph associated to the McKeithan scheme is strongly connected and with
      deficiency zero. In each strictly positive compatibility class there exists in this
      domain an unique globally asymptotically stable equilibrium.
      \label{thm:anderson}
    \end{theorem}

    \subsection{Application to stabilization and observer design for the McKeithan scheme}

    This stability result has consequences to control and observation properties of the 
    network, see \cite{sontag2001,chavez2002b}, that we recall briefly.
    The dynamics \eqref{eq:dyn} can be converted into a control system of the form
  \begin{equation}
    \begin{aligned}
      &\dot c(t) = f(x(t)) + u(t)
      \\
      &y(t) = h(c(t)),
    \end{aligned}
    \label{eq:control-sys}
  \end{equation}
  where $u(\cdot)$ is a feeding control and $h$ is a polynomial observation function.

  Asymptotic stability of equilibrium in each strictly positive compatibility class will 
  allow to get stabilization result for a single equilibrium.
  Moreover it leads to design under a mild assumption (detectability) a simple observer.
  We refer to \cite{sontag2001} and \cite{chavez2002a} for the detailed presentation of those
  results, the geometric construction being clear.

  \section{The optimal control problem and Pontryagin Maximum Principle}
    \label{sec3}

  \subsection{Statement and notation for the optimal control problem}
  \label{sec3p1}
  The system is written as 
    $\frac{\dd c}{\dd t} = f(c,T)$
    (see \eqref{eq:dyn})
    and controlling the temperature leads to $T\in [T_m,T_M]$. 
    In the sequel, we shall use the terminology {\it direct} for the corresponding 
    control problem. In practice, thermodynamics has to be used to model the heat exchanges,
    in relation with the heat produced by the reactions \cite{toth2018} or the heat exchange
    device used in the experiments, depending upon the technical achievements. To avoid
    this part of the study and without losing any mathematically generality, we shall use
    $\dot v$ as the control variable setting $\dot v=u$, where $v\coloneqq k_i(T)$
    for some reaction $i$, where $k_i(T)=A_i e^{-E_i/(RT)}$ (see \eqref{eq:constant-reaction}).
      This leads to deal with the so--called {\it indirect} control system:
      \[
        \dot q(t) = F(q(t)) + u(t)\, G(q(t)),
      \]
      where $q=(c,v) \in \mathbb{R}^n$ is the extended state variable, 
      $F=f(c,v)\, \frac{\partial}{\partial q}, G=\frac{\partial}{\partial v}$, 
      $u_-\le u \le u_+$, where $[u_-,u_+]$ can be normalized to $[-1,+1]$. Note that
      the bounds $v\in [v_m,v_M]$ will not be taken into account in our study.
      The map $v\mapsto \dot v$ is the standard {\it Goh transformation} in optimal control,
      see \cite{bonnard2003}.

      The optimal control problem of physical interest is the problem of maximizing the 
      production of one species and using a proper variable labeling, the optimal problem
      is therefore of the {\it Mayer type}:
      \[
        \dot q = F(q) + u\, G(q),\quad \underset{|u|\le 1}{\max}\ q_1(t_f),
      \]
      where $t_f$ is the time duration of the batch and $q_1$ is the desired product.
      
      Note it will be show (thanks to the Maximum Principle) that a ``dual'' formulation
      is 
      \[
        \underset{u(\cdot)}{\min} \ t_f, \quad c_1(t_f) = d,
      \]
      where $d>0$ is the desired amount of the species $X_1$ during the batch.

      \subsection{Maximum Principle \cite{pontryagin1962}}
      \subsubsection{Notations and concepts}
      Consider a general control system of the form
      \[
        \dot q = X(q,u), \quad q\in \mathbb{R}^n,
      \]
      where $X$ is a real analytic $(C^\omega)$ and the control is 
      $u:[0,t_f(u)] \mapsto [-1,1]$.
      The set of admissible controls ${\cal U}$ is the set of bounded measurable 
      mappings. 
      If $q(0)=q_0$ (initial state), we denote by $q(\cdot,q_0,u)$ (in short $q(\cdot)$)
      the solution starting from $q_0$. Fixing $t_f$, the {\it accessibility set}
      in time $t_f$ is the set 
      $A(q_0,t_f) = \underset{u(\cdot)\in {\cal U}}{\cup} q(t_f,q_0,u)$.
      The {\it extremity mapping} (in time $t_f$) is the map:
      $E^{q_0,t_f}: u(\cdot) \mapsto q(t_f,q_0,u)$ defined on a domain 
      of ${\cal U}$; the set ${\cal U}$ is endowed with the L$^\infty$-norm
      topology.

      \subsubsection{Pontryagin Maximum Principle (PMP)}

\noindent
      {\bf Statement in the time minimal case with} $\bm{q(t_f)\in N\coloneqq}$ 
        {\bf smooth terminal manifold.}

        \noindent
        {\it Notation: } $H(q,p,u) = p\cdot X(q,u)$ denotes the {\it pseudo--Hamiltonian}
        (Hamiltonian lift of the vector field $X$), $p$ is the adjoint vector in 
        $\mathbb{R}^n\setminus \{0\}$ and $\cdot$ is the scalar product. 
        We denote by $M(q,p) = \underset{|u|\le 1}{\max}\ H(q,p,u)$.

\noindent
        {\it Statement of the PMP:} If $(q^*,u^*)$ is an optimal control--trajectory 
        pair on $[0,t_f^*]$ then there exists $p^*\in \mathbb{R}^n\setminus \{0\}$
        such that a.e.:
        \begin{equation}
          \begin{array}{ll}
            &\dot q^*(t)=\frac{\partial H}{\partial p}(q^*(t),p^*(t),u^*(t)),
            \\
            &\dot p^*(t) = -\frac{\partial H}{\partial q}(q^*(t),p^*(t),u^*(t)),
            \\
            &H(q^*(t),p^*(t),u^*(t))=M(q^*(t),p^*(t)).
          \end{array}
        \label{eq:ham-pmp}
        \end{equation}
        Moreover $M(q^*(t),p^*(t))$ is a positive constant and $p^*$ satisfies
        the {\it transversality condition}
\begin{equation}
p^*(t_f) \perp T_{q^*(t_f)} N.
\label{eq:transversality}
\end{equation}   

\noindent
{\it Statement in the Mayer case:} As before, except that the transversality condition is 
replaced by 
\[
  p^*(t_f)= \frac{\partial \varphi}{\partial q}(x^*(t_f)),
\]
where $\varphi$ is the Mayer cost function to minimize.

\begin{definition}
  An extremal $(q,p,u)$ is a solution of \eqref{eq:ham-pmp} on $[0,t_f]$. It is called a 
  BC--extremal is $q(0)=q_0$ and $p(t_f)$ satisfies the transversality condition.
  An extremal control is called regular if $|u(t)|\le 1$ a.e. and singular if
  $\frac{\partial H}{\partial u}=0$ everywhere.
  An extremal is said exceptional if $M=0$. 
  A regular extremal control is called "bang--bang" if $u(\cdot)$ is piecewise constant
  on $[0,t_f]$ (i.e. the number of switches is finite).
\end{definition}

\subsubsection{Computations of singular extremals}

\noindent
{\bf First case.} Consider the case $\dot q = X(q,u)$, where $H=p\cdot X(q,u)$ and the condition
 $\frac{\partial H}{\partial u}=0$ is satisfied. 
Denote by $z(\cdot)=(q(\cdot),p(\cdot))$ the reference extremal. 
From the maximization condition,
the {\it Legendre--Clebsch condition} $\frac{\partial^2 H}{\partial u^2}\le 0$ has to be fulfilled.
If this inequality is strict, one can use the implicit function theorem to compute the 
singular control as the dynamic feedback: $z\mapsto u_s(z)$ and plugging such $u_s$
into $H(q,p,u)$ leads to define the {\it true} (or maximized) Hamiltonian.

\noindent
{\bf Second case.} Let $\dot q = F(q) + u\, G(q)$. One introduces the following notations.
If $X$ and $Y$ are two real analytic vector fields, the {\it Lie bracket} is defined by
$[X,Y](q)= \frac{\partial X}{\partial q}(q)Y(q)-\frac{\partial Y}{\partial q}(q)X(q)$. 
The extremal lift of $X$ is $H_X(z)=p\cdot X(q), \; z=(q,p)$ and the 
{\it Poisson bracket} is defined by 
$\{H_X,H_Y \}(z)=p\cdot [X,Y](q)$.

\noindent
{\it Computation of singular extremals:}

The condition 
$\frac{\partial}{\partial u} \frac{\dd}{\dd t^2} \frac{\partial H}{\partial u}
=\{ \{ H_G,H_F \}, H_G \} \geq 0$ (resp. $>0$) 
is called the (resp. strict) {\it generalized Legendre--Clebsch condition}.

Recall the following.

\begin{proposition}[\cite{krener1977}]
  The generalized Legendre--Clebsch condition is a necessary optimality condition for the 
  time minimal control problem with fixed extremities.
 \end{proposition} 
  To compute the singular extremal, we differentiate twice $t\mapsto H_G(z(t))$
  and we get 
  \begin{equation}
  \begin{aligned}
&H_G(z)=\{ H_G,H_F \}(z)=0,
\\
&\{ \{ H_G,H_F \} , H_F \}(z)+u\, \{ \{ H_G,H_F \} , H_G \}(z)=0. 
  \end{aligned}
  \label{eq:cond_order_2}
  \end{equation}

Assume (in relation with the Legendre--Clebsch condition) ${\{\{H_G,H_F\},H_G\}\neq 0}$.
The corresponding extremal is called of {\it order 2} and the singular control $u$ is
computed as $u_s(z)$, using relation \eqref{eq:cond_order_2}.
Plugging such $u_s(z)$ into $H(q,p,u)$ leads to define the {\it true}
singular Hamiltonian, denoted by $H_s(z)$. One has:

\begin{proposition}
  Singular extremals of order $2$ are the solutions of :
\begin{equation}
\frac{\dd q}{\dd t}=\frac{\partial H_s}{\partial p}(z), 
\; \frac{\dd p}{\dd t}=-\frac{\partial H_s}{\partial q}(z)
\label{eq:singular_solution_Hs}
\end{equation}
\end{proposition}
with the constraints 
\[\{ H_G,H_F \}(z) =  H_G(z) = 0.\]
Moreover, in order to be admissible, the singular control is given by
\begin{equation}
u_s(z)= -\frac{\{ \{ H_G,H_F \} , H_F \}(z)}{\{ \{ H_G,H_F \} , H_G \}(z)}
\label{eq:control_u_s}
\end{equation} 
and has to satisfy the admissibility constraint $|u_s(z)|\le 1$.

\begin{definition}
  Let $z=(q,p)$ be a singular extremal of order $2$  and $M=H_F=h$ the constant value of 
  the Hamiltonian.
  The extremal is called {\it exceptional} if $h=0$. If $h>0$, 
  the extremal is called 
  {\it hyperbolic} (resp. {\it elliptic}) if 
  $\{\{H_G,H_F\},H_G\}>0$ (resp. $<0$).
  \label{def:class-sing}
\end{definition}

\subsubsection{The case of chemical networks}

Recall that for our network $\dot c = f(c,v),\ v=k_i$ and it is extended into 
${\dot q = F(q) + u\, G(q)}$ with $F=f(c,v)\frac{\partial}{\partial c}$ and
$G=\frac{\partial}{\partial v}$. 
Denote $\tilde H = p_c\cdot f(c,v)$ and $H=p\cdot (F+u\, G)$ 
the respective Hamiltonian lifts, with $p=(p_c,p_v)$. One has the following
relation between the corresponding singular extremals.

\begin{lemma}
  The pair $(q,p)$ is solution of 
  \[
    \dot q =  \frac{\partial H}{\partial p},\ 
    \dot p = -\frac{\partial H}{\partial q},\
    \frac{\partial H}{\partial u}=0
  \]
  if and only if $p_v=0$ and $(c,p_c,p_v)$ is solution of 
  \[
    \dot{c}=\frac{\partial \tilde{H}}{\partial p_c}, \; 
    \dot{p}_c=\frac{\partial \tilde{H}}{\partial \tilde{c}}, \; 
    \frac{\partial \tilde{H}}{\partial v}=0,
  \] 
  and moreover the following relations are satisfied   
\begin{equation}
\frac{\dd}{\dd t}\frac{\partial H}{\partial u}=\{H_G,H_F \}=-\frac{\partial \tilde{H}}{\partial v}
\label{eq:H_tilde_v}
\end{equation}
\begin{equation}
\frac{\partial }{\partial u}\frac{\dd^2}{\dd t^2}\frac{\partial H}{\partial u}=\{ \{H_G,H_F \}, H_G \}=-\frac{\partial^2 \tilde{H}}{\partial v^2}.
\label{eq:H_tilde_v2}
\end{equation}
\end{lemma}

In particular this gives the correspondence between both singular Hamiltonian 
and the respective Legendre--Clebsch and generalized Legendre--Clebsch conditions.

\subsubsection{The case $n=3$}
We have $q=(x,y,v)$ and introduce the following determinants:
\begin{align*}
&D=\det(G,[G,F],[[G,F],G]),
\\
& D'=\det(G,[G,F],[[G,F],F]),
\\
&D''=\det(G,[G,F],F).
\end{align*}
Using $p\cdot G = p\cdot [G,F]= p\cdot \left(\left[\left[ G,F \right],F\right]
+ u \, \left[\left[ G,F \right],G\right]\right)=0$
and eliminating $p$ leads to the following.

\begin{proposition}
  If $n=3$, the singular control is given by the feedback ${u_s(q)=-D'(q)/D(q)}$
  and singular trajectories are defined by the vector field 
  $X_s(q) = F(q) - \frac{D'(q)}{D(q)}\, G(q)$. 
  Singular trajectories are hyperbolic if $DD''>0$, elliptic if $DD''<0$ and 
  exceptional if $D''=0$.
\end{proposition}

\noindent
{\bf Optimality status: the case $\bm{n=3}$.}
We use \cite{bb1993} to describe the optimality status of singular trajectories.
The system is $\dot q = F+u\, G$ and we relax the bound $|u|\le 1$, assuming 
$u\in \mathbb{R}$ so that singular arcs are admissible. We assume the following:
\begin{description}
  \item[($H_0$)] $F,G$ are linearly independent, 
  \item[($H_1$)] $G,[G,F]$ are linearly independent. 
\end{description}
One picks a smooth singular arc $z(\cdot) = (q(\cdot),p(\cdot))$ defined on $[0,t_f]$
so that $p$ is unique up to a non zero multiplicative scalar 
and $q(\cdot)$ is a one-to-one immersion. We have:
\begin{proposition}
\label{prop:classification-singular}
  There exists a $C^0$--neighborhood $U$ of the reference singular arc 
  ${t\mapsto q(t),\ t\in [0,t_f]}$ so that $q(\cdot)$ is time minimal (resp. maximal)
  if $q(\cdot)$ is hyperbolic (resp. elliptic) up to the first conjugate time $t_{1c}$ 
  with respect to all trajectories with the same extremities contained in ${\cal U}$.
  In the exceptional case, the reference singular arc is time minimal (and time 
  maximal).
\end{proposition}

\noindent
{\bf Algorithm to compute the first conjugate point in the hyperbolic--elliptic
case.}
Let $V(t),\ t\in[0,t_f]$ be the solution of the {\it (variational) Jacobi equation}:
\[
  \dot{\delta q}(t) = \frac{\partial X}{\partial q}(z(t))\, \delta q(t)
\]
with initial condition $V(0)=G(q(0))$, where $t\rightarrow q(t)$ is the reference
singular arc. 
The first conjugate point $t_{1c}$ is the first $t>0$ such that  $V(t)$ is collinear to
$G(q(t))$.

See \cite[p. 123]{bonnard2003} for a proof and the geometric interpretation.

\subsection{Small time classification of regular extremals}
\label{sec:class-extremals}

In this section, we recall the seminal results coming from singularity theory due to
\cite{ekeland1977,kupka1987} to analyze the small time extremal curves near the switching surface.

\begin{definition}
  We denote by $\sigma_+$ (resp. $\sigma_-$) a bang arc with constant control $u=1$ (resp.
  $u=-1$) and $\sigma_s$ an admissible singular arc. We denote by $\sigma_1 \sigma_2$
  an arc $\sigma_1$ followed by $\sigma_2$.
  The surface $\Sigma: H_G(z)=0$ is called the switching surface and let $\Sigma'\subset \Sigma$
  given by $H_G(z)=\{H_G,H_F\}(z)=0$.
  Let $z(\cdot)=(q(\cdot),p(\cdot))$ be a reference curve on $[0,t_f]$.
  We note $t\mapsto \Phi(t)\coloneqq H_G(z(t))$ the switching function, which codes the switching 
  times.
\end{definition}
 
 Deriving twice with respect to time the switching function, one gets
\begin{equation}
  \dot{\Phi}(t)=\{H_G,H_F\}(z(t))\,,
\label{eq:Phi_dot}
\end{equation}   
\begin{equation}
  \ddot{\Phi}(t)=\{\{H_G,H_F\},H_F\}(z(t))+u(t)\, \{\{H_G,H_F\},H_G\}(z(t)).
\label{eq:Phi_ddot}
\end{equation}   

From this, we derive.
\begin{lemma}
\label{lem:classification-ordinary}
  Assume that $t$ is an ordinary switching time, that is $\Phi(t)=0$ and 
  ${\dot \Phi(t)\neq 0}$. Then near $z(t)$ every extremal projects onto $\sigma_+\sigma_-$
  if $\dot \Phi(t)<0$ and $\sigma_-\sigma_+$ if $\dot \Phi(t)>0$.
\end{lemma}

The situation is more complex for higher contact with $\Sigma$.
\begin{definition}
\label{def:fold_case}
  The fold case is characterized by $\Phi(t)=\dot \Phi(t)=0$ and 
  $\ddot \Phi(t)\neq 0$ (replacing $u$ by $\pm 1$ in \eqref{eq:Phi_ddot}).
  Hence $z(t)\in \Sigma'$.
  Assume that locally $\Sigma'$ is a regular surface of codimension two. We
  have three cases:
  \begin{itemize}
    \item 
      {\it Parabolic case:} 
      $\ddot{\Phi}_{+}(t) \ddot{\Phi}_{-}(t) > 0$.
    \item 
      {\it Hyperbolic case:} 
      $\ddot{\Phi}_{+}(t)>0$ and $\ddot{\Phi}_{-}(t) < 0$.
    \item 
      {\it Elliptic case:}
      $\ddot{\Phi}_{+}(t)<0$ and $\ddot{\Phi}_{-}(t) > 0$. 
  \end{itemize}
\end{definition}
Denote by $u_s(\cdot)$ the singular control defined by \eqref{eq:Phi_ddot} as
$\ddot \Phi(t)=0$. In the hyperbolic case, through $z(t)$, assuming
$\left\{ \left\{ H_G,H_F \right\},H_G \right\}(z(t))\neq 0$, there exists a singular arc, which 
is strictly admissible that is $|u_s(t)|<1$. This arc is hyperbolic
if  $\left\{ \left\{ H_G,H_F \right\},H_G \right\}(z(t))>0$ (strict Legendre--Clebsch condition)
and elliptic if this quantity is $<0$. In the parabolic case, it can be absent or not admissible
that is $|u_s(t)|>1$.

One has from \cite{kupka1987},
\begin{figure}
    \centering
    \def\svgwidth{0.44\textwidth}
%% Accompanies image file 'fold_case_elliptic.pdf' (pdf, eps, ps)
\begingroup%
  \makeatletter%
  \providecommand\color[2][]{%
    \errmessage{(Inkscape) Color is used for the text in Inkscape, but the package 'color.sty' is not loaded}%
    \renewcommand\color[2][]{}%
  }%
  \providecommand\transparent[1]{%
    \errmessage{(Inkscape) Transparency is used (non-zero) for the text in Inkscape, but the package 'transparent.sty' is not loaded}%
    \renewcommand\transparent[1]{}%
  }%
  \providecommand\rotatebox[2]{#2}%
  \newcommand*\fsize{\dimexpr\f@size pt\relax}%
  \newcommand*\lineheight[1]{\fontsize{\fsize}{#1\fsize}\selectfont}%
  \ifx\svgwidth\undefined%
    \setlength{\unitlength}{252.66439603bp}%
    \ifx\svgscale\undefined%
      \relax%
    \else%
      \setlength{\unitlength}{\unitlength * \real{\svgscale}}%
    \fi%
  \else%
    \setlength{\unitlength}{\svgwidth}%
  \fi%
  \global\let\svgwidth\undefined%
  \global\let\svgscale\undefined%
  \makeatother%
  \begin{picture}(1,0.42128937)%
    \lineheight{1}%
    \setlength\tabcolsep{0pt}%
    \put(0,0){\includegraphics[width=\unitlength,page=1]{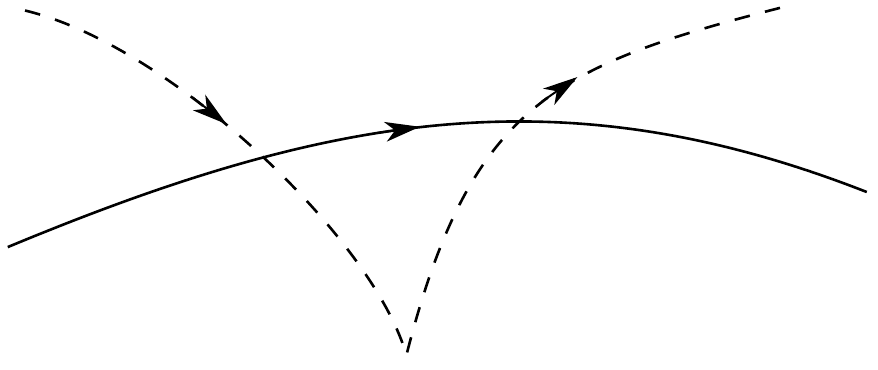}}%
    \put(0.22038985,0.33433285){\color[rgb]{0,0,0}\makebox(0,0)[lt]{\lineheight{1.25}\smash{\begin{tabular}[t]{l}$\sigma_+$\end{tabular}}}}%
    \put(0.58887272,0.347729){\color[rgb]{0,0,0}\makebox(0,0)[lt]{\lineheight{1.25}\smash{\begin{tabular}[t]{l}$\sigma_-$\end{tabular}}}}%
    \put(0.40746358,0.20829872){\color[rgb]{0,0,0}\makebox(0,0)[lt]{\lineheight{1.25}\smash{\begin{tabular}[t]{l}$\sigma_s$\end{tabular}}}}%
  \end{picture}%
\endgroup%

    \caption{Fold case in the elliptic 2D--situation and the antiturnpike phenomenon.}
    \label{fig:Fold_case_turnpike}
\end{figure}
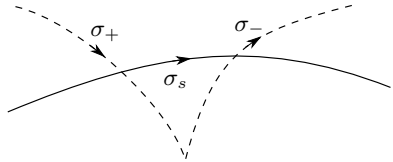

\subsection{Global case} 
\begin{figure}
    \centering
    \def\svgwidth{0.42\textwidth}
%% Accompanies image file 'global_case.pdf' (pdf, eps, ps)
\begingroup%
  \makeatletter%
  \providecommand\color[2][]{%
    \errmessage{(Inkscape) Color is used for the text in Inkscape, but the package 'color.sty' is not loaded}%
    \renewcommand\color[2][]{}%
  }%
  \providecommand\transparent[1]{%
    \errmessage{(Inkscape) Transparency is used (non-zero) for the text in Inkscape, but the package 'transparent.sty' is not loaded}%
    \renewcommand\transparent[1]{}%
  }%
  \providecommand\rotatebox[2]{#2}%
  \newcommand*\fsize{\dimexpr\f@size pt\relax}%
  \newcommand*\lineheight[1]{\fontsize{\fsize}{#1\fsize}\selectfont}%
  \ifx\svgwidth\undefined%
    \setlength{\unitlength}{352.50193919bp}%
    \ifx\svgscale\undefined%
      \relax%
    \else%
      \setlength{\unitlength}{\unitlength * \real{\svgscale}}%
    \fi%
  \else%
    \setlength{\unitlength}{\svgwidth}%
  \fi%
  \global\let\svgwidth\undefined%
  \global\let\svgscale\undefined%
  \makeatother%
  \begin{picture}(1,0.50972748)%
    \lineheight{1}%
    \setlength\tabcolsep{0pt}%
    \put(0,0){\includegraphics[width=\unitlength,page=1]{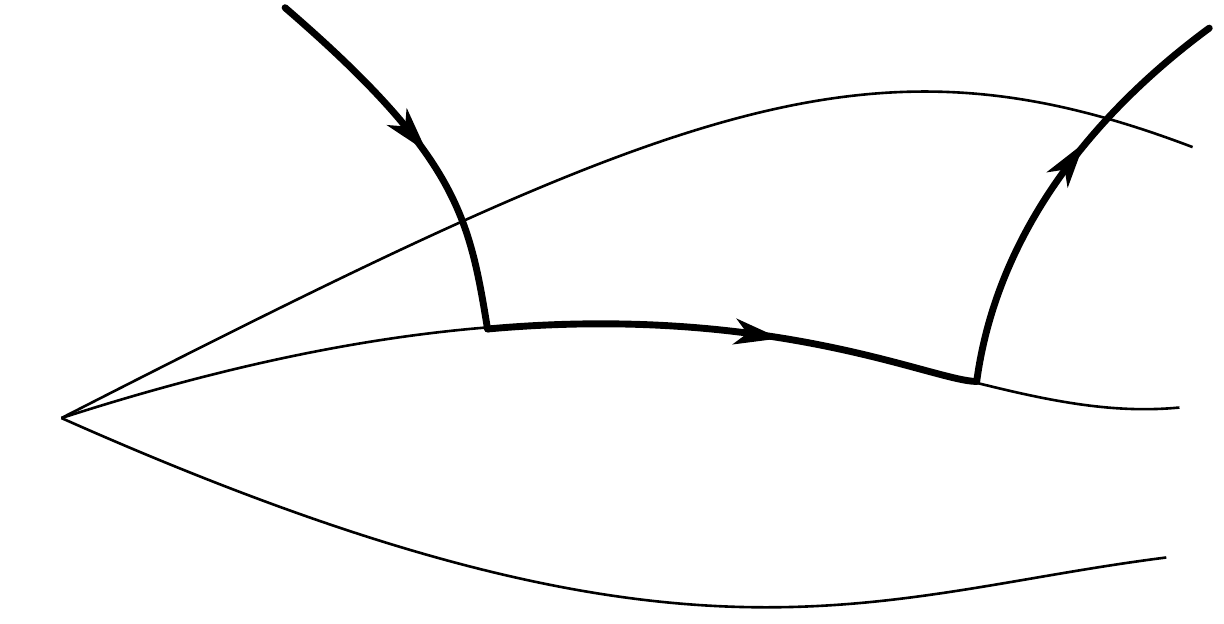}}%
    \put(-0.08129701,0.13618665){\color[rgb]{0,0,0}\makebox(0,0)[lt]{\lineheight{1.25}\smash{\begin{tabular}[t]{l}$q(0)$\end{tabular}}}}%
    \put(0.49564859,0.18819752){\color[rgb]{0,0,0}\makebox(0,0)[lt]{\lineheight{1.25}\smash{\begin{tabular}[t]{l}$\sigma_s$\end{tabular}}}}%
    \put(0.30567469,0.43947126){\color[rgb]{0,0,0}\makebox(0,0)[lt]{\lineheight{1.25}\smash{\begin{tabular}[t]{l}$\sigma_+$\end{tabular}}}}%
    \put(0.857614,0.30587856){\color[rgb]{0,0,0}\makebox(0,0)[lt]{\lineheight{1.25}\smash{\begin{tabular}[t]{l}$\sigma_-$\end{tabular}}}}%
    \put(0,0){\includegraphics[width=\unitlength,page=2]{global_case.pdf}}%
    \put(0.40273828,0.12050123){\color[rgb]{0,0,0}\rotatebox{-10.282489}{\makebox(0,0)[lt]{\lineheight{1.25}\smash{\begin{tabular}[t]{l}singular arcs\end{tabular}}}}}%
    \put(0,0){\includegraphics[width=\unitlength,page=3]{global_case.pdf}}%
  \end{picture}%
\endgroup%

    \def\svgwidth{0.42\textwidth}
%% Accompanies image file 'conj-c-space.pdf' (pdf, eps, ps)
\begingroup%
  \makeatletter%
  \providecommand\color[2][]{%
    \errmessage{(Inkscape) Color is used for the text in Inkscape, but the package 'color.sty' is not loaded}%
    \renewcommand\color[2][]{}%
  }%
  \providecommand\transparent[1]{%
    \errmessage{(Inkscape) Transparency is used (non-zero) for the text in Inkscape, but the package 'transparent.sty' is not loaded}%
    \renewcommand\transparent[1]{}%
  }%
  \providecommand\rotatebox[2]{#2}%
  \newcommand*\fsize{\dimexpr\f@size pt\relax}%
  \newcommand*\lineheight[1]{\fontsize{\fsize}{#1\fsize}\selectfont}%
  \ifx\svgwidth\undefined%
    \setlength{\unitlength}{334.54139049bp}%
    \ifx\svgscale\undefined%
      \relax%
    \else%
      \setlength{\unitlength}{\unitlength * \real{\svgscale}}%
    \fi%
  \else%
    \setlength{\unitlength}{\svgwidth}%
  \fi%
  \global\let\svgwidth\undefined%
  \global\let\svgscale\undefined%
  \makeatother%
  \begin{picture}(1,0.4037383)%
    \lineheight{1}%
    \setlength\tabcolsep{0pt}%
    \put(0,0){\includegraphics[width=\unitlength,page=1]{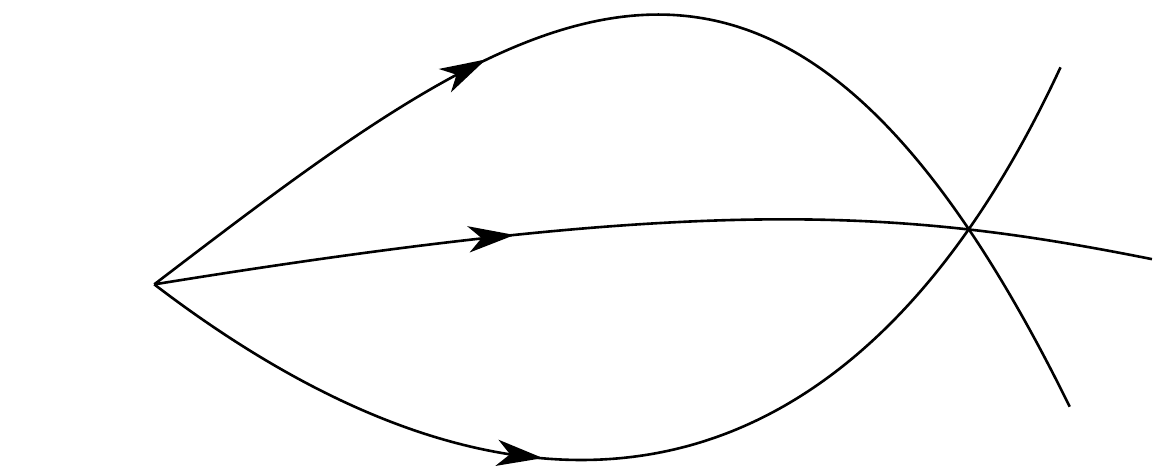}}%
    \put(0.44071287,0.16173671){\color[rgb]{0,0,0}\makebox(0,0)[lt]{\lineheight{1.25}\smash{\begin{tabular}[t]{l}$\sigma_s$\end{tabular}}}}%
    \put(0.72827268,0.22950217){\color[rgb]{0,0,0}\makebox(0,0)[lt]{\lineheight{1.25}\smash{\begin{tabular}[t]{l}$t_{1c}$\end{tabular}}}}%
    \put(0.02056693,0.14409576){\color[rgb]{0,0,0}\makebox(0,0)[lt]{\lineheight{1.25}\smash{\begin{tabular}[t]{l}$c(0)$\end{tabular}}}}%
    \put(0,0){\includegraphics[width=\unitlength,page=2]{conj-c-space.pdf}}%
  \end{picture}%
\endgroup%

  \caption{ {\it(left)} Time minimal policy in the $q$--space. {\it (right)} Conjugate 
point in the $c$--space.}
    \label{fig:global-case}
\end{figure}

It is based on \cite{bb1993} and presented here 
for $n=3$. It uses proposition \ref{prop:classification-singular}. 
Under our assumptions $(H_0)-(H_1)$,
and if moreover the singular control is {\it strictly admissible} on $[0,t_f]$,
that is $|u_s(t)|<1$, the reference singular arc can be immersed in a 
$C^0$--domain {\it up to the first conjugate time} $t_{1c}$ and where the time
minimal policy is of the form $\sigma_\pm \sigma_s \sigma_\pm$, see Fig. \ref{fig:global-case}
for the interpretation of the first conjugate time $t_{1c}$ for the problem in the
$q$--space with $q=(c,v)$, $\dot v = u$.

\subsection{The curse of the non strict Legendre--Clebsch condition}

More complicated and challenging situation is to analyze situations, where the strict 
Legendre--Clebsch condition is not satisfied that is
$\left\{ \left\{ H_G,H_F \right\},H_G \right\}(z(t))$ vanishes at some times.
Due to the complexity and in relation with our application, we shall assume that
$n=3$.

In this case, one has $D=0$ since 
$\left\{ \left\{ H_G,H_F \right\},H_G \right\} = p \cdot \left[ \left[ G,F \right],G \right]$
and recall that $D=\det(G,[G,F],\left[ \left[ G,F \right],G \right])$.
The singular flow is defined by the vector field 
$X_s(q)=F(q) - \frac{D'(q)}{D(q)}\, G(q)$. Using a time reparameterization, it is related
to the one--dimensional foliation ${\cal F}_s: D(q)F(q)-D'(q)G(q)$.
Complexity of the analysis is due to {\it non isolated singularities} contained in the
set $D'(q)=D(q)=0$.

In relation with our optimal problem with terminal manifold  $N$ of codimension one given
by $x=d$, our singularity resolution will be concerned by analyzing arcs
initiating from the set ${\cal S}:n \cdot \left[ G,F \right](q)=0$, where $n=(1,0,0)$ is the 
normal vector to $N$. Singularities can be roughly classified into two types: local in 
relation with singularities of ${\cal S}$ and propagated along the singular flow, and Lagrangian
singularities in relation with the concept of conjugate--focal points.
This will be developed in the next section.

\section{Geometric techniques to analyze the 2d--McKeithan Network}
\label{sec4}
\subsection{Time minimal syntheses near the terminal manifold}

\subsubsection{Notations and definitions}

The system is written $\dot q = F(q) + u\, G(q)$, $|u|\le 1$, with
$q=(c,v)$, $G=\frac{\partial}{\partial v}$.
The terminal manifold $N$ is given by $c_1=d$. The problem is to determine 
{\it small time} synthesis near a given point $q_0\in N$, which can be identified to $0$.
More precisely one wants to classify the syntheses under generic assumptions near the
terminal manifold in relation with the Lie algebraic structure of $\left\{ F,G \right\}$
at $q_0$. Our approach developed in our series of articles \cite{bonnard1995,bonnard1997,bonnard1998,launay1997} is to use the
construction of {\it semi--normal form} for the action of the pseudo--group ${\cal G}$
of local diffeomorphisms and feedback transformation $u\rightarrow-u$ (so that
$\sigma_+$ and $\sigma_-$ can be exchanged). 
Additionally recall that the pseudo--group ${\cal G}_f$ formed by local diffeomorphims
and feedback actions of the form $u=\alpha(q) + \beta(q)v$ leaves the singular flow
invariant \cite{bb1991}. These groups act on the jet space of $F$ at zero, $G$ being identified
to $\frac{\partial}{\partial v}$ and the terminal manifold $N$ to $c_1=d$. 
Note that the problem is {\it flat} that is $G$ is tangent to $N$.
For the action of the pseudo--group $G$ on the jet spaces of $F$ we refer to \cite{martinet1980}
(we shall work in the $C^k$ category, where $k\ge 1$ is not precised) and the 
semi--normal form is related to a semi--algebraic stratification on the jet spaces of
$(F,G,N)$.
One has $N: c_1=d$ and the initial state is such that $c_1(0)<d$. Denote in general 
$N^\perp = \{(q,p); \ p\cdot v=0, \forall v \in T_qN\}$ and let $n$ be the outward 
normal to $N$ so that $n=(1,0,\ldots,0)$ and $N^\perp = (q,n(q))$.
Let $z=(q,p)$ be a BC--extremal on $[t_f,0]$, $t_f<0$ and $z(0)\in N^\perp$ 
(this convention is used since we integrate {\it backwards} from the terminal manifold).
The set of minimizing switching points can be stratified into stratum of 
{\it first kind} if the optimal curves are tangent
and {\it second kind} if they are transverse.
The {\it splitting locus} $L$ is the set of points, where the optimal control is not unique
and the {\it cut locus} $C$ is the closure of the set of points, where the optimal
trajectories loses its optimality, see \cite{boltyanskii1966,brunovsky1980} for the introduction of those
concepts in the frame of semi--analytic geometry.

We introduce the following triplet $(F,G,N)$ in the $C^\omega$--category, 
with $N=\{f(q)=0\}$, where $f$ is a (local) submersion. Fixing $q_0\in \mathbb{R}^n$,
we denote by $j_kF(q_0),j_kG(q_0),j_kf(q_0)$ the respective $k$--jets, that is the
Taylor expansions at order $k$. We say that $(F,G,f)$ has at $q_0$ a singularity
of {\it codimension} $i$ if $(j_kF(q_0),j_kG(q_0),j_kf(q_0)) \in \Sigma_i$,
where $\Sigma_i$ is a semi--algebraic submanifold of codimension $i$ in the 
jets space.

Taking $q_0\in N$ with a singularity of codimension $i$, an {\it unfolding} is a
$C^0$ change of coordinates near $q_0$ such that small time minimal synthesis
is described by a system $\dot {\tilde q} = F(\tilde q, \lambda) + u\, G(\tilde q,\lambda)$,
$|u|\le 1$, $\tilde q\in \mathbb{R}^{n-m}$ and $\lambda$ is a (vector) parameter.

\subsubsection{Local syntheses: general tools and classification}
We shall present the main steps to compute the time minimal syntheses, restricting
our study to the $3$-dimensional case, but it can be clearly extended to the 
$n$--space with the concept of codimension \cite{launay1997}.

The system is written: $\dot q = F(q) + u\, G(q)$, $|u|\le 1$ and 
let $q=(x,y,z)$ be the coordinates, $G$ being identified to $\frac{\partial}{\partial z}$
and $G$ being tangent to $N$, which can be identified to $x=0$ and $n=(1,0,0)$
is the outward normal to $N$. Recall that (generic) singular trajectories are given
by $\dot q = F(q) + u_s(q) G(q)$, $u_s(q)=-D'(q)/D(q)$ and they can be classified into:
hyperbolic, elliptic, exceptional, see section \ref{sec3}.

The first step is to stratify the terminal manifold into: 
\begin{itemize}
\item
${\cal S}:$ {\it singular locus}
defined by $\{q\in N,\ n\cdot [G,F](q)=0\}$,
\item
${\cal E}:$ {\it exceptional locus}
defined by $\{q\in N,\ n\cdot F(q)=0\}$.
\end{itemize}

Note that since the problem is flat, $n\cdot G(q)=0$ if $q\in N$ so that 
$N^\perp \subset \Sigma$, where $\Sigma$ is the switching surface.

\noindent
{\bf Generic case.} The first case is when $F(q_0)$ and $\left[ G,F \right](q_0)$ are 
independent and not tangent to $N$ and the synthesis follows from Lemma 
\ref{lem:classification-ordinary}.
It is given by $\sigma_+$ if $n\cdot \left[ G,F \right](q_0)<0$ and $\sigma_-$ if
$n\cdot [G,F](q_0)>0$ since the final point is a {\it virtual} switching point.

\noindent
{\bf Generic hyperbolic singular case and the concept of focal points.}
Now we present the basis of our analysis, the so--called hyperbolic situation. In this case,
based on \cite{bonnard1997}, a semi--algebraic normal form is constructed to obtain the
corresponding local syntheses and it is extended (in the jet space) {\it along 
the reference singular arc} in order to obtain the concept of {\it focal point} 
based on the
notions of turnpike and conjugate points presented in Section \ref{sec3}.

\noindent
{\it Semi--normal form:} one has $q_0=0$ and we make the following assumptions:
\begin{itemize}
  \item 
    The tangent space to $N$ is $G(0)$, $\left[ G,F \right](0)$. 
  \item The set of points $\delta$, 
  where $\left[ G,F \right]$ is tangent to $N$ is a simple
    curve passing through $0$ and transverse to $G$.
  \item 
    $D(0)$ and $D''(0)$ are non zero. 
\end{itemize}
With those assumptions, through $0$, there exists 
    a simple BC--singular extremal $\sigma_s$ transverse to $N$. One can choose local 
    coordinates so that $G$ is identified to $\frac{\partial}{\partial z}$, $\delta$
    to the axis $(Oy)$ can be identified to $t\mapsto (t,0,0)$, whose image is the
    $(Ox)$--axis. Note that $N$ is identified to $x=0$.

The semi--normal form is constructed in a tubular neighborhood of the reference singular
curve $\sigma_s$ and the vector field $F$ is developed in the jet space with respect to
$(y,z)$ since $\sigma_s:t\mapsto (t,0,0)$ is the x--coordinate. One has the 
following semi-normal form
\begin{align*}
&\dot{x}=1+a(x) z^2+2 b(x) y z +c(x) y^2+R_1 
\\
&\dot{y}=d(x) y+ e(0) z+ R_2 
\\
&\dot{z}=(u-{u_s}_{\mid\sigma}(x))+ f(x) y+ g(0) z + R_3,
\end{align*}  
where $a(x)\neq 0, \; e(0) \neq 0$ and $R_1$ 
(resp. $R_2, \; R_3$) are terms of order $\geq 3$ 
(resp. $\geq 2$) in $(y,z)$ and $u_s$ is the singular control.

Furthermore we make the following assumption:
\begin{itemize}
  \item The reference arc is {\it hyperbolic} on $[t_f,0]$ so that $a(x)<0$.
\end{itemize}

\begin{figure}
    \centering
    \def\svgwidth{0.3\textwidth}
%% Accompanies image file 'synthesis-bb.pdf' (pdf, eps, ps)
\begingroup%
  \makeatletter%
  \providecommand\color[2][]{%
    \errmessage{(Inkscape) Color is used for the text in Inkscape, but the package 'color.sty' is not loaded}%
    \renewcommand\color[2][]{}%
  }%
  \providecommand\transparent[1]{%
    \errmessage{(Inkscape) Transparency is used (non-zero) for the text in Inkscape, but the package 'transparent.sty' is not loaded}%
    \renewcommand\transparent[1]{}%
  }%
  \providecommand\rotatebox[2]{#2}%
  \newcommand*\fsize{\dimexpr\f@size pt\relax}%
  \newcommand*\lineheight[1]{\fontsize{\fsize}{#1\fsize}\selectfont}%
  \ifx\svgwidth\undefined%
    \setlength{\unitlength}{370.4911648bp}%
    \ifx\svgscale\undefined%
      \relax%
    \else%
      \setlength{\unitlength}{\unitlength * \real{\svgscale}}%
    \fi%
  \else%
    \setlength{\unitlength}{\svgwidth}%
  \fi%
  \global\let\svgwidth\undefined%
  \global\let\svgscale\undefined%
  \makeatother%
  \begin{picture}(1,0.56030034)%
    \lineheight{1}%
    \setlength\tabcolsep{0pt}%
    \put(0,0){\includegraphics[width=\unitlength,page=1]{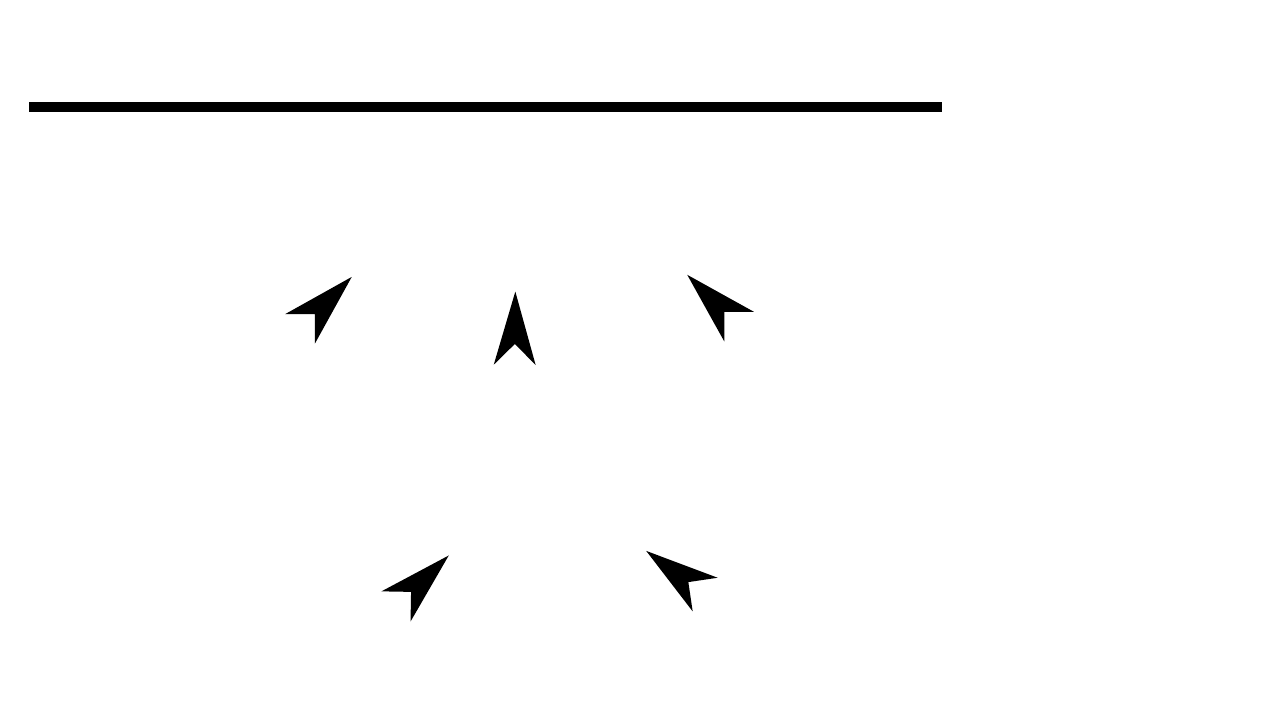}}%
    \put(0.37218912,0.49942278){\color[rgb]{0,0,0}\makebox(0,0)[lt]{\lineheight{1.25}\smash{\begin{tabular}[t]{l}$0$\end{tabular}}}}%
    \put(0.13602255,0.37335348){\color[rgb]{0,0,0}\makebox(0,0)[lt]{\lineheight{1.25}\smash{\begin{tabular}[t]{l}$\sigma_+$\end{tabular}}}}%
    \put(0.58561793,0.33354967){\color[rgb]{0,0,0}\makebox(0,0)[lt]{\lineheight{1.25}\smash{\begin{tabular}[t]{l}$\sigma_-$\end{tabular}}}}%
    \put(0.42855028,0.25766996){\color[rgb]{0,0,0}\makebox(0,0)[lt]{\lineheight{1.25}\smash{\begin{tabular}[t]{l}$\sigma_s$\end{tabular}}}}%
    \put(0.74825139,0.46197732){\color[rgb]{0,0,0}\makebox(0,0)[lt]{\lineheight{1.25}\smash{\begin{tabular}[t]{l}$x=0$\end{tabular}}}}%
    \put(0.75484273,0.17133988){\color[rgb]{0,0,0}\makebox(0,0)[lt]{\lineheight{1.25}\smash{\begin{tabular}[t]{l}$x<0$\end{tabular}}}}%
    \put(0,0){\includegraphics[width=\unitlength,page=2]{synthesis-bb.pdf}}%
  \end{picture}%
\endgroup%

    \caption{}
    \label{fig:synthesis-bb}
\end{figure}

According to \cite{bonnard1997} the time minimal synthesis near $q_0$ is of the form
$\sigma_+ \sigma_s \sigma_-$, see Fig.\ref{fig:synthesis-bb} and is described by the
unfolding 
\begin{align*}
  &\dot x = 1 - a(0)\, z^2
  \\
  &\dot z = u - u_s(0)
  \\
  &\dot y = 0
\end{align*}
in a small neighborhood of $q_0=0$.

\paragraph*{Local syntheses and the semi-bridge phenomenon.}

We use \cite{bonnard1997} and the \texttt{Mathematica}'s program given 
in Appendix \ref{prog:strata}.
One fixes a point $q_0 \in {\cal S}$ and we make the explicit computations of the
switching and cut loci near $q_0$ using the jet expansion of $F$ at $q_0$,
$G$ and $N$ being normalized respectively to $\frac{\partial}{\partial v}$
and $x=0$.
The computations are in the semi-algebraic category, easily  
implementable using symbolic algorithms and the cut locus and switching loci
are stratifiable.

{\it Switching locus.} We denote by $K$ the set of ordinary switching points
for BC-extremals, $K_+$ being a switching $\sigma_- \sigma_+$ and $K_-$ being
a switching $\sigma_+ \sigma_-$, while $W,W_+,W_-$ corresponds respectively to
switching points of optimal extremals. 
More precisely, near $N$ the stratification of $W$ is $W=W_1 \cup W_2$ 
where $W_1$ is the first kind stratum 
composed by the hyperbolic singular arcs
and
$W_2=W_s\cup W_+ \cup W_-$ is the second kind strata, where
$W_\varepsilon,\ \varepsilon\in \{-1,1\}$ 
is composed by the ordinary switching points of the policy 
$\sigma_{-\varepsilon} \sigma_\varepsilon$ and $W_s$ is the first switching
point of the Bang-Bang-Singular's policy.

\begin{figure}[ht]
  \centering
    \def\svgwidth{0.25\textwidth}
%% Accompanies image file 'ka.pdf' (pdf, eps, ps)
\begingroup%
  \makeatletter%
  \providecommand\color[2][]{%
    \errmessage{(Inkscape) Color is used for the text in Inkscape, but the package 'color.sty' is not loaded}%
    \renewcommand\color[2][]{}%
  }%
  \providecommand\transparent[1]{%
    \errmessage{(Inkscape) Transparency is used (non-zero) for the text in Inkscape, but the package 'transparent.sty' is not loaded}%
    \renewcommand\transparent[1]{}%
  }%
  \providecommand\rotatebox[2]{#2}%
  \newcommand*\fsize{\dimexpr\f@size pt\relax}%
  \newcommand*\lineheight[1]{\fontsize{\fsize}{#1\fsize}\selectfont}%
  \ifx\svgwidth\undefined%
    \setlength{\unitlength}{180.26785663bp}%
    \ifx\svgscale\undefined%
      \relax%
    \else%
      \setlength{\unitlength}{\unitlength * \real{\svgscale}}%
    \fi%
  \else%
    \setlength{\unitlength}{\svgwidth}%
  \fi%
  \global\let\svgwidth\undefined%
  \global\let\svgscale\undefined%
  \makeatother%
  \begin{picture}(1,0.6225855)%
    \lineheight{1}%
    \setlength\tabcolsep{0pt}%
    \put(0,0){\includegraphics[width=\unitlength,page=1]{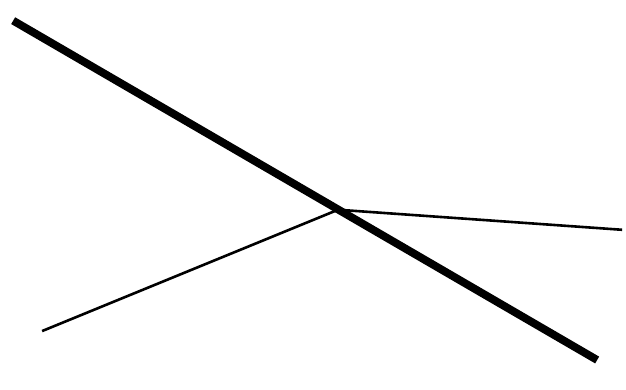}}%
    \put(0.20335048,0.50761377){\color[rgb]{0,0,0}\makebox(0,0)[lt]{\lineheight{1.25}\smash{\begin{tabular}[t]{l}$K$\end{tabular}}}}%
    \put(0.11600998,0.03495681){\color[rgb]{0,0,0}\makebox(0,0)[lt]{\lineheight{1.25}\smash{\begin{tabular}[t]{l}$\sigma_+$\end{tabular}}}}%
    \put(0.72069209,0.33178695){\color[rgb]{0,0,0}\makebox(0,0)[lt]{\lineheight{1.25}\smash{\begin{tabular}[t]{l}$\sigma_-$\end{tabular}}}}%
    \put(0,0){\includegraphics[width=\unitlength,page=2]{ka.pdf}}%
  \end{picture}%
\endgroup%

\hspace{1cm}
       \def\svgwidth{0.25\textwidth}
%% Accompanies image file 'kb.pdf' (pdf, eps, ps)
\begingroup%
  \makeatletter%
  \providecommand\color[2][]{%
    \errmessage{(Inkscape) Color is used for the text in Inkscape, but the package 'color.sty' is not loaded}%
    \renewcommand\color[2][]{}%
  }%
  \providecommand\transparent[1]{%
    \errmessage{(Inkscape) Transparency is used (non-zero) for the text in Inkscape, but the package 'transparent.sty' is not loaded}%
    \renewcommand\transparent[1]{}%
  }%
  \providecommand\rotatebox[2]{#2}%
  \newcommand*\fsize{\dimexpr\f@size pt\relax}%
  \newcommand*\lineheight[1]{\fontsize{\fsize}{#1\fsize}\selectfont}%
  \ifx\svgwidth\undefined%
    \setlength{\unitlength}{180.26785663bp}%
    \ifx\svgscale\undefined%
      \relax%
    \else%
      \setlength{\unitlength}{\unitlength * \real{\svgscale}}%
    \fi%
  \else%
    \setlength{\unitlength}{\svgwidth}%
  \fi%
  \global\let\svgwidth\undefined%
  \global\let\svgscale\undefined%
  \makeatother%
  \begin{picture}(1,0.6225855)%
    \lineheight{1}%
    \setlength\tabcolsep{0pt}%
    \put(0,0){\includegraphics[width=\unitlength,page=1]{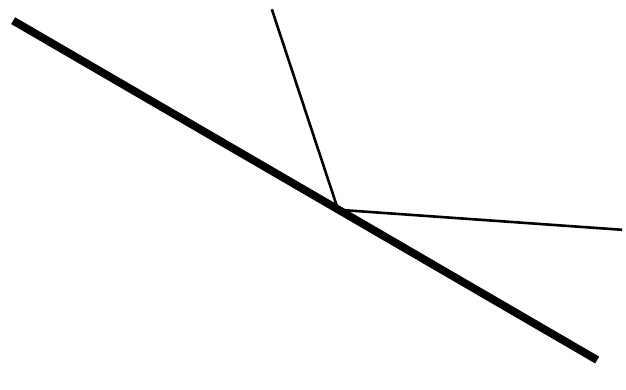}}%
    \put(0.20335048,0.50761377){\color[rgb]{0,0,0}\makebox(0,0)[lt]{\lineheight{1.25}\smash{\begin{tabular}[t]{l}$K$\end{tabular}}}}%
    \put(0.46667861,0.55204445){\color[rgb]{0,0,0}\makebox(0,0)[lt]{\lineheight{1.25}\smash{\begin{tabular}[t]{l}$\sigma_-$\end{tabular}}}}%
    \put(0.72069209,0.33178695){\color[rgb]{0,0,0}\makebox(0,0)[lt]{\lineheight{1.25}\smash{\begin{tabular}[t]{l}$\sigma_+$\end{tabular}}}}%
    \put(0,0){\includegraphics[width=\unitlength,page=2]{kb.pdf}}%
  \end{picture}%
\endgroup%

  \caption{ {\it (left)} Crossing. {\it (right)} Reflecting.
\label{fig:crossing-k}}
\end{figure}

From \cite{Agrachev1998} we have two situations describe by Fig.\ref{fig:crossing-k}.
The reflecting case corresponds to a non optimal situation and has to be rejected to
determine $W$.	

Using the theory and computations of \cite{bonnard1997} only the stratification of
$W$ localized near $q_0$ can be computed.

{\it Singular locus.} Near $q_0$, one determines the singular locus $\Gamma_s$
restricting to admissible singular trajectories, which are optimal. 
In other words, the test is :
hyperbolicity and strict admissibility $|u_s|<1$.

{\it Cut locus.} The cut locus $C$ is the set of points where optimality is lost.
One part of the cut locus is formed by the splitting locus $L$, where two minimizers 
intersects. It can be stratified into strata corresponding between intersections $\sigma_+$, $\sigma_-$
intersections between $\sigma_+ \sigma_-$ and $\sigma_-$ \dots Again such intersections
are described in \cite{bonnard1997}.

\paragraph*{The local syntheses.} It is based on the classification of section \ref{sec:class-extremals}.

{\bf Generic case.} Since $G$ is tangent to $N$, every final point is a virtual switching point.
Denoting $\Phi(t)$ the switching function on $[t_f,0]$, with $\Phi(0)=n\cdot G(q_0)$ one has:
if $\dot \Phi(0)\coloneqq n\cdot[G,F](q_0)<0$ (resp. $>0$) the terminating arc is $\sigma_-$
(resp. $\sigma_+$).

\begin{figure}[ht]
    \begin{minipage}{0.49\textwidth}
    \centering
        \def\svgwidth{0.5\textwidth}
  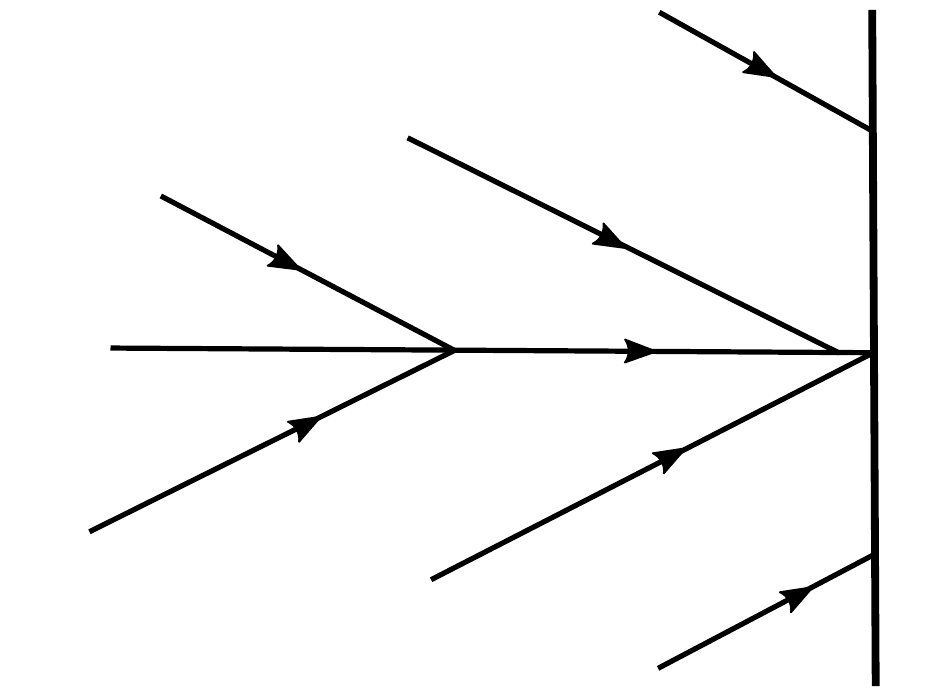
   \\
    Hyperbolic
\end{minipage}
\begin{minipage}{0.49\textwidth}
    \centering
    \def\svgwidth{0.37\textwidth}
%% Accompanies image file 'parabolic_elliptic.pdf' (pdf, eps, ps)
\begingroup%
  \makeatletter%
  \providecommand\color[2][]{%
    \errmessage{(Inkscape) Color is used for the text in Inkscape, but the package 'color.sty' is not loaded}%
    \renewcommand\color[2][]{}%
  }%
  \providecommand\transparent[1]{%
    \errmessage{(Inkscape) Transparency is used (non-zero) for the text in Inkscape, but the package 'transparent.sty' is not loaded}%
    \renewcommand\transparent[1]{}%
  }%
  \providecommand\rotatebox[2]{#2}%
  \newcommand*\fsize{\dimexpr\f@size pt\relax}%
  \newcommand*\lineheight[1]{\fontsize{\fsize}{#1\fsize}\selectfont}%
  \ifx\svgwidth\undefined%
    \setlength{\unitlength}{225.53569163bp}%
    \ifx\svgscale\undefined%
      \relax%
    \else%
      \setlength{\unitlength}{\unitlength * \real{\svgscale}}%
    \fi%
  \else%
    \setlength{\unitlength}{\svgwidth}%
  \fi%
  \global\let\svgwidth\undefined%
  \global\let\svgscale\undefined%
  \makeatother%
  \begin{picture}(1,1.04988124)%
    \lineheight{1}%
    \setlength\tabcolsep{0pt}%
    \put(0,0){\includegraphics[width=\unitlength,page=1]{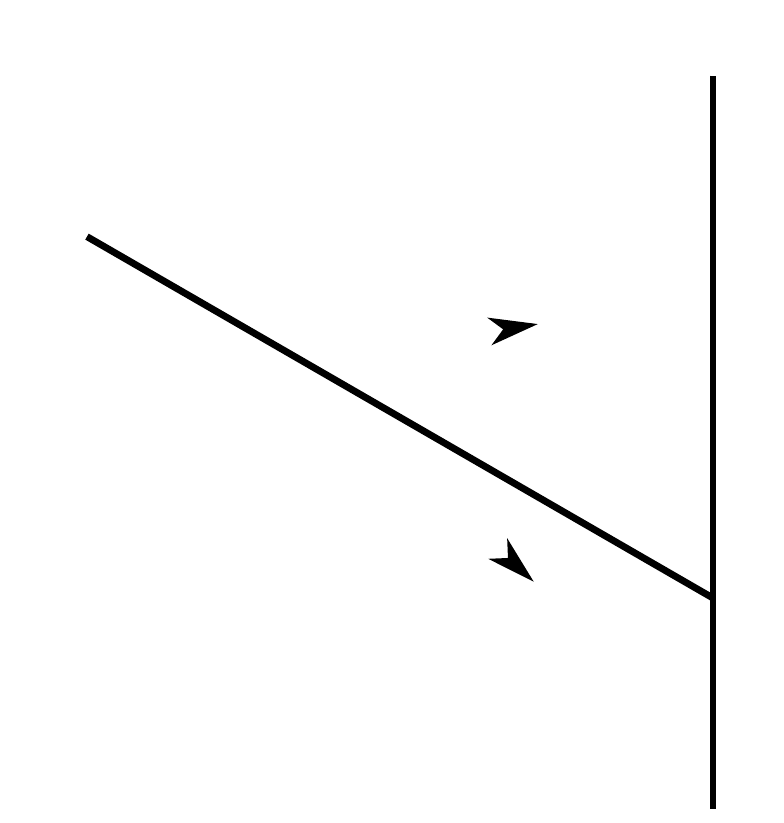}}%
    \put(0.00475049,0.76128252){\color[rgb]{0,0,0}\makebox(0,0)[lt]{\lineheight{1.25}\smash{\begin{tabular}[t]{l}$C$\end{tabular}}}}%
    \put(0.84898275,0.97191462){\color[rgb]{0,0,0}\makebox(0,0)[lt]{\lineheight{1.25}\smash{\begin{tabular}[t]{l}$N$\end{tabular}}}}%
    \put(0.50255595,0.25829295){\color[rgb]{0,0,0}\makebox(0,0)[lt]{\lineheight{1.25}\smash{\begin{tabular}[t]{l}$\sigma_-$\end{tabular}}}}%
    \put(0.50824127,0.68886366){\color[rgb]{0,0,0}\makebox(0,0)[lt]{\lineheight{1.25}\smash{\begin{tabular}[t]{l}$\sigma_+$\end{tabular}}}}%
    \put(0,0){\includegraphics[width=\unitlength,page=2]{parabolic_elliptic.pdf}}%
  \end{picture}%
\endgroup%

   \\
  \quad Elliptic
\end{minipage}
    \caption{Fold case.
  \label{fig:hyperbolic_case}}
\end{figure}

\begin{figure}[ht]
  \centering
    \begin{minipage}{0.49\textwidth}
    \centering
    \def\svgwidth{0.35\textwidth}
   
   \hspace{0.5cm}
       \def\svgwidth{0.35\textwidth}
%% Accompanies image file 'parabolic_elliptic-b.pdf' (pdf, eps, ps)
\begingroup%
  \makeatletter%
  \providecommand\color[2][]{%
    \errmessage{(Inkscape) Color is used for the text in Inkscape, but the package 'color.sty' is not loaded}%
    \renewcommand\color[2][]{}%
  }%
  \providecommand\transparent[1]{%
    \errmessage{(Inkscape) Transparency is used (non-zero) for the text in Inkscape, but the package 'transparent.sty' is not loaded}%
    \renewcommand\transparent[1]{}%
  }%
  \providecommand\rotatebox[2]{#2}%
  \newcommand*\fsize{\dimexpr\f@size pt\relax}%
  \newcommand*\lineheight[1]{\fontsize{\fsize}{#1\fsize}\selectfont}%
  \ifx\svgwidth\undefined%
    \setlength{\unitlength}{225.53569163bp}%
    \ifx\svgscale\undefined%
      \relax%
    \else%
      \setlength{\unitlength}{\unitlength * \real{\svgscale}}%
    \fi%
  \else%
    \setlength{\unitlength}{\svgwidth}%
  \fi%
  \global\let\svgwidth\undefined%
  \global\let\svgscale\undefined%
  \makeatother%
  \begin{picture}(1,1.04988124)%
    \lineheight{1}%
    \setlength\tabcolsep{0pt}%
    \put(0,0){\includegraphics[width=\unitlength,page=1]{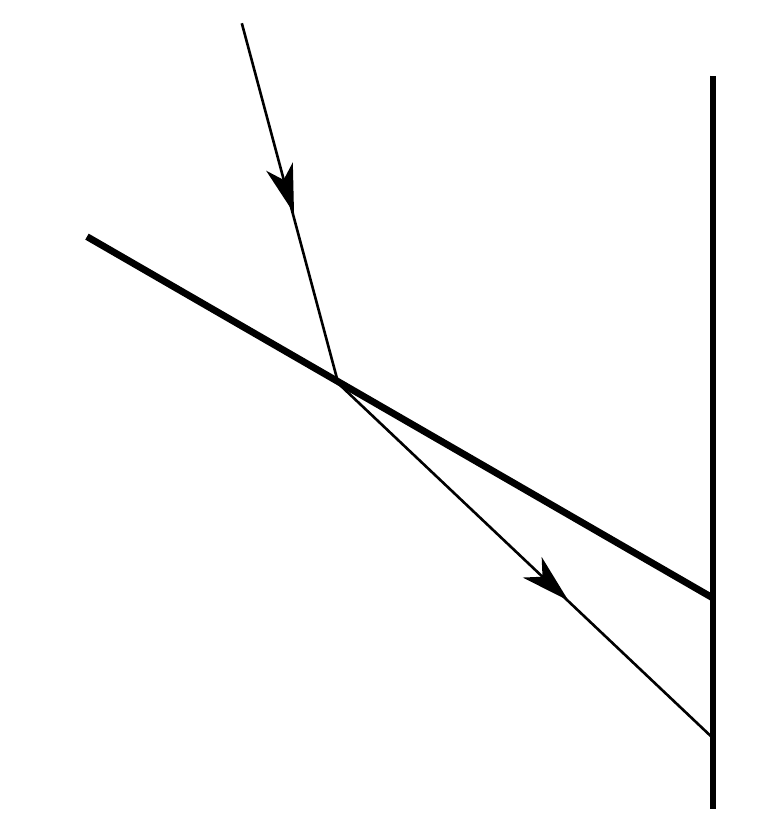}}%
    \put(0.00475049,0.76128252){\color[rgb]{0,0,0}\makebox(0,0)[lt]{\lineheight{1.25}\smash{\begin{tabular}[t]{l}$W_-$\end{tabular}}}}%
    \put(0.84898275,0.97191462){\color[rgb]{0,0,0}\makebox(0,0)[lt]{\lineheight{1.25}\smash{\begin{tabular}[t]{l}$N$\end{tabular}}}}%
    \put(0.34898256,0.96953923){\color[rgb]{0,0,0}\makebox(0,0)[lt]{\lineheight{1.25}\smash{\begin{tabular}[t]{l}$\sigma_+$\end{tabular}}}}%
    \put(0.45325743,0.32702138){\color[rgb]{0,0,0}\makebox(0,0)[lt]{\lineheight{1.25}\smash{\begin{tabular}[t]{l}$\sigma_-$\end{tabular}}}}%
    \put(0,0){\includegraphics[width=\unitlength,page=2]{parabolic_elliptic-b.pdf}}%
    \put(0.64536194,0.71538273){\color[rgb]{0,0,0}\makebox(0,0)[lt]{\lineheight{1.25}\smash{\begin{tabular}[t]{l}$\sigma_+$\end{tabular}}}}%
    \put(0,0){\includegraphics[width=\unitlength,page=3]{parabolic_elliptic-b.pdf}}%
  \end{picture}%
\endgroup%

   \\
   Non admissible elliptic case
\end{minipage}
\begin{minipage}{0.49\textwidth}
  \centering
    \def\svgwidth{0.35\textwidth}
%% Accompanies image file 'parabolic_hyperbolic.pdf' (pdf, eps, ps)
\begingroup%
  \makeatletter%
  \providecommand\color[2][]{%
    \errmessage{(Inkscape) Color is used for the text in Inkscape, but the package 'color.sty' is not loaded}%
    \renewcommand\color[2][]{}%
  }%
  \providecommand\transparent[1]{%
    \errmessage{(Inkscape) Transparency is used (non-zero) for the text in Inkscape, but the package 'transparent.sty' is not loaded}%
    \renewcommand\transparent[1]{}%
  }%
  \providecommand\rotatebox[2]{#2}%
  \newcommand*\fsize{\dimexpr\f@size pt\relax}%
  \newcommand*\lineheight[1]{\fontsize{\fsize}{#1\fsize}\selectfont}%
  \ifx\svgwidth\undefined%
    \setlength{\unitlength}{225.53569163bp}%
    \ifx\svgscale\undefined%
      \relax%
    \else%
      \setlength{\unitlength}{\unitlength * \real{\svgscale}}%
    \fi%
  \else%
    \setlength{\unitlength}{\svgwidth}%
  \fi%
  \global\let\svgwidth\undefined%
  \global\let\svgscale\undefined%
  \makeatother%
  \begin{picture}(1,1.04988124)%
    \lineheight{1}%
    \setlength\tabcolsep{0pt}%
    \put(0,0){\includegraphics[width=\unitlength,page=1]{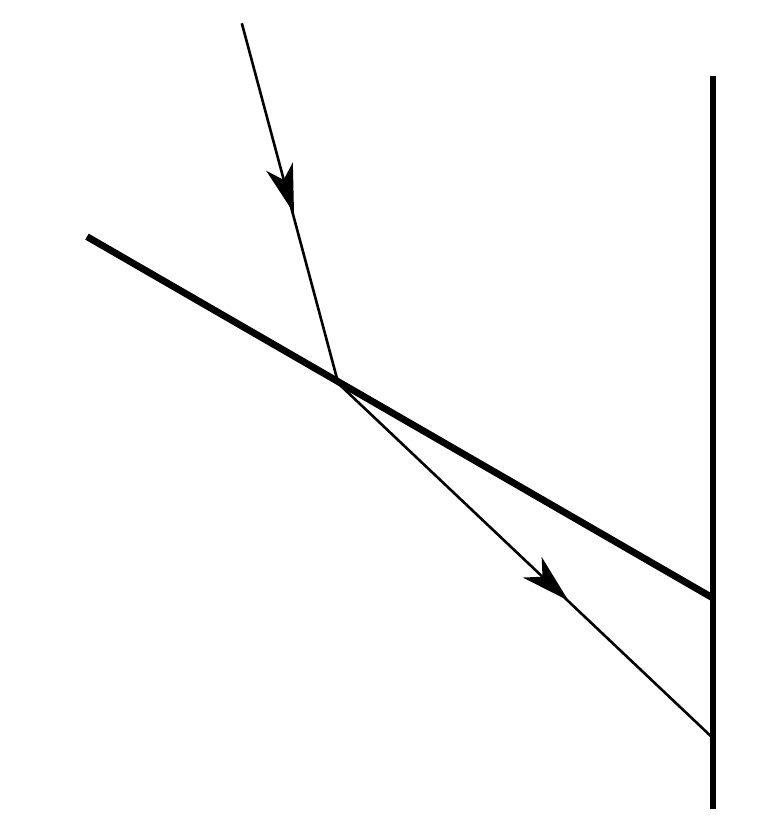}}%
    \put(0.03800468,0.77553427){\color[rgb]{0,0,0}\makebox(0,0)[lt]{\lineheight{1.25}\smash{\begin{tabular}[t]{l}$W_+$\end{tabular}}}}%
    \put(0.84898275,0.97191462){\color[rgb]{0,0,0}\makebox(0,0)[lt]{\lineheight{1.25}\smash{\begin{tabular}[t]{l}$N$\end{tabular}}}}%
    \put(0.37154788,0.91965799){\color[rgb]{0,0,0}\makebox(0,0)[lt]{\lineheight{1.25}\smash{\begin{tabular}[t]{l}$\sigma_-$\end{tabular}}}}%
    \put(0.45325743,0.32702138){\color[rgb]{0,0,0}\makebox(0,0)[lt]{\lineheight{1.25}\smash{\begin{tabular}[t]{l}$\sigma_+$\end{tabular}}}}%
    \put(0,0){\includegraphics[width=\unitlength,page=2]{parabolic_hyperbolic.pdf}}%
  \end{picture}%
\endgroup%

   \\
   Non admissible hyperbolic case
 \end{minipage}
  \caption{Parabolic point with a non admissible singular arc.
\label{fig:parabolic_point}}
\end{figure}

{\bf Codimension one case.}
We have different cases corresponding to the fold case, since in the parabolic case one must 
distinguish between a parabolic point corresponding either to a non admissible 
hyperbolic or elliptic arc.
The cases are represented on Fig.\ref{fig:hyperbolic_case}-\ref{fig:parabolic_point}.

The syntheses can be represented by foliations by 2d-planes. Note that the role of
$\sigma_+$ and $\sigma_-$ can be interchanged since the semi-normal forms are computed 
using the transformation $u\mapsto -u$.

The ${\cal C}^0$-unfoldings are:
\begin{itemize}
  \item Hyperbolic and Parabolic:
    \begin{align*}
    &\dot x = 1 + a\, y^2
    \\
    &\dot y = u - u_s(0), \qquad |u_s(0)|<1,
    \end{align*}
    and $a<0$ (hyperbolic) or $a>0$ (elliptic).	
%  \item Elliptic:
%    \begin{align*}
%    &\dot x = 1 + a\, y^2
%    \\
%    &\dot y = u - u_s(0), \qquad |u_s(0)|<1
%    \end{align*}
  \item Parabolic:
    \begin{align*}
    &\dot x = 1 + a\, y^2
    \\
    &\dot y = u - u_s(0), \qquad |u_s(0)|>1,
    \end{align*}
    and $a\neq 0$.
\end{itemize}

\paragraph*{Focal points.}
The synthesis in the hyperbolic case can be extended to a tubular 
neighborhood of $\sigma_s$ for $t\in [t_f,0]$
introducing the concept of focal point as follows.
By assumptions, $W\coloneqq \delta'(0)$ belongs to $\mathrm{span} \{G(0),[G,F](0)\}$.
Let $\lambda_1, \lambda_2$ be two scalars such that 
$W=\lambda_1 G(0) + \lambda_2 [F,G](0)$ and by assumption $\lambda_1\neq 0$.
Denoting by $X_s(q)\coloneqq F(q) - \frac{D'(q)}{D(q)} G(q)$, we introduce the variational
equation 
\[
  \dot{\delta q}(t) = \frac{\partial X_s}{\partial q}(\sigma_s(t))\, \delta q(t)
\]
and let $W(\cdot)$ be a solution on $[t_f,0]$ such that $W(0)=W$.
\begin{definition}
  Let $t_{1f}$ be the first time in $[t_f,0]$ such that:
  ${\det(W(t),G(\sigma_s(t)),F(\sigma_s(t)))=0}$. 
  Then $t_{1f}$ is called the first focal point along $\sigma_s$.
\end{definition}
As for the fixed end point problem we have the following:
\begin{proposition}
Provided $t\in ]t_{1f},0]$, then in a tubular neighborhood of $\sigma_s$ the set
$S=\{\exp(t X_s(\sigma))\}$ is a smooth surface and the synthesis is given by 
Fig.\ref{fig:synthesis-bb}.
\end{proposition}
\noindent
{\it Algorithm:} Besides the definition, which leads to compute focal point using 
a singular value decomposition, an equivalent computation is to determine $t=t_{1f}$
such that $W(t)$ becomes collinear to $G(\sigma_s(t))$.

{\bf Two codimension-two cases.}

The situation is more intricate. It is analyzed using the semi-normal form constructed in
\cite{bonnard1997}.
A  model  is
\begin{equation}
\begin{aligned}
&\dot x = 1 + a\, z^2+\alpha_1\, x y^2+\alpha_2\, y z^2+\alpha_3\, x z^2
    \\
&\dot y = b\, z
    \\
&\dot z = c\, z-u_s(0)-u_{sx}\, x-u_{sy}\, y,
\end{aligned}
\label{eq:semi-nf}
\end{equation}
where $u_s$ is the singular control,
$u_s(0) \in \{1,3\}$,
$u_{sx}=\frac{\partial u_s}{\partial x}(0)$,
$u_{sy}=\frac{\partial u_s}{\partial y}(0)>0$ and $a,b,c\neq 0$. 
The switching function $\Phi(t)\coloneqq p_3(t)$ is developed at order $3$ and factorized as
$tP(t)$, where $P$ is polynomial of order two with two roots $t_1,t_2$, which determine the
switching points.

If $a>0$, we are in the elliptic situation and if $a<0$ in the hyperbolic situation and the reference
singular arc is identified to $\sigma_s:t\mapsto (t,0,0)$ being not admissible.
One has $|u_s(0)|>1$ and we can assume $u_s(0)>1$.

Figures presented in this section are obtained with the \texttt{Mathematica}'s program
of Appendix \ref{prog:strata}, we use the following values:
$\alpha_1=\alpha_2=\alpha_3=u_{sy}=c=b=1$.\\

{\bf The case $\bm{a>0}$ and $\bm{u_s(0)=3}$.} We have the following three situations
\begin{figure}[H]
    \centering
\begin{tabular}{ccc}    
    \def\svgwidth{0.25\textwidth}
   
   &
         \def\svgwidth{0.3\textwidth}
   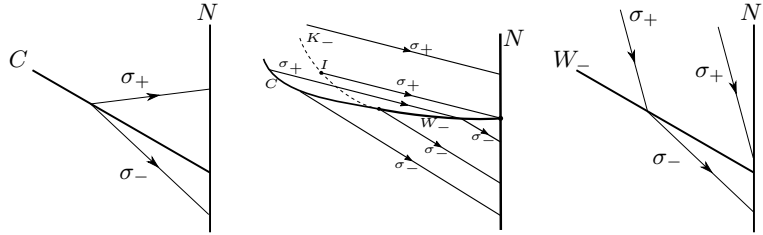 \hspace*{0.2cm}
   &
       \def\svgwidth{0.25\textwidth}
   
   \end{tabular}
   \caption{\label{fig:par-ell}
Bifurcation of the cut $(a>0,\, u_s(0)=3)$.
 The cut locus $C$ splits into $C_1\cup C_{12}$ where $C_1: \sigma_-\cap\sigma_+$
 and $C_{12}: \sigma_+\sigma_-\cap \sigma_-$.}
\end{figure}

\begin{figure}[H]
  \centering
\begin{tabular}{c}
\includegraphics[scale=0.7]{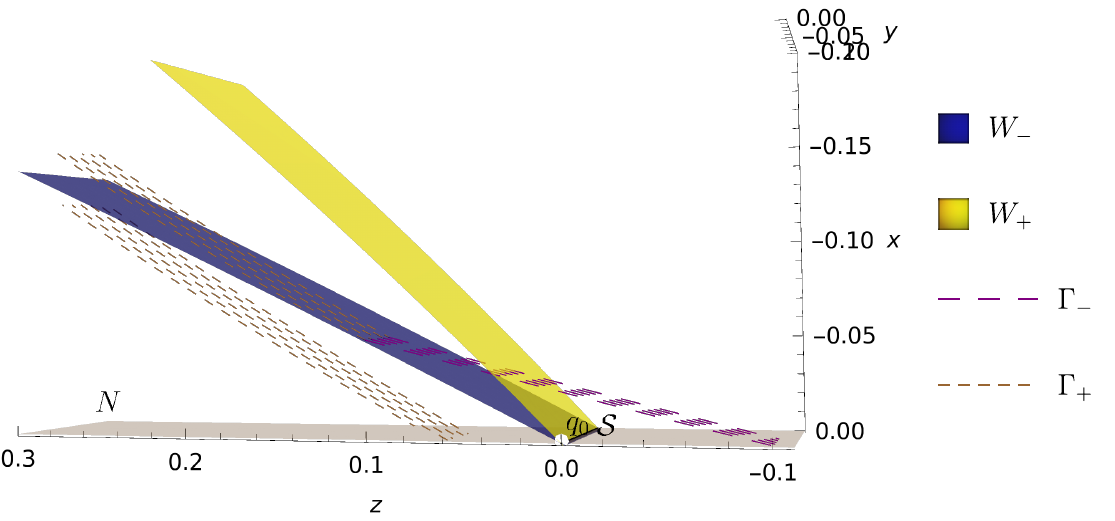} 
\\
$y<0$
\end{tabular}
\caption{
  A situation in the elliptic parabolic case.
  The trajectories of $\Gamma_-$ reflect on $W_-$
leading to a cut locus between $W_-$ and $W_+$.\label{fig:cut-elliptic}}
\end{figure}

Note that the cut appears when the trajectories reflect on
the switching surface, see Fig.\ref{fig:cut-elliptic}.

%The transition between the two cases of Fig.\ref{fig:par-ell} is for the specific value $u_s_0=3$
%in the normal form.
%In particular it is a situation where the cut has  represented in the Fig.\ref{fig:cut-decolle}.

\begin{figure}
\includegraphics[scale=0.47]{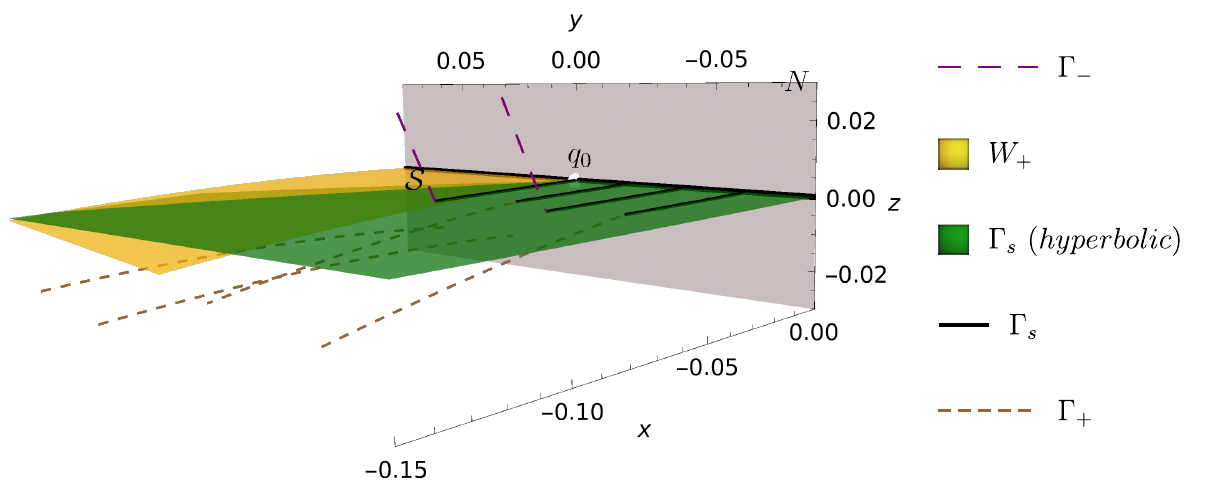} 
\includegraphics[scale=0.47]{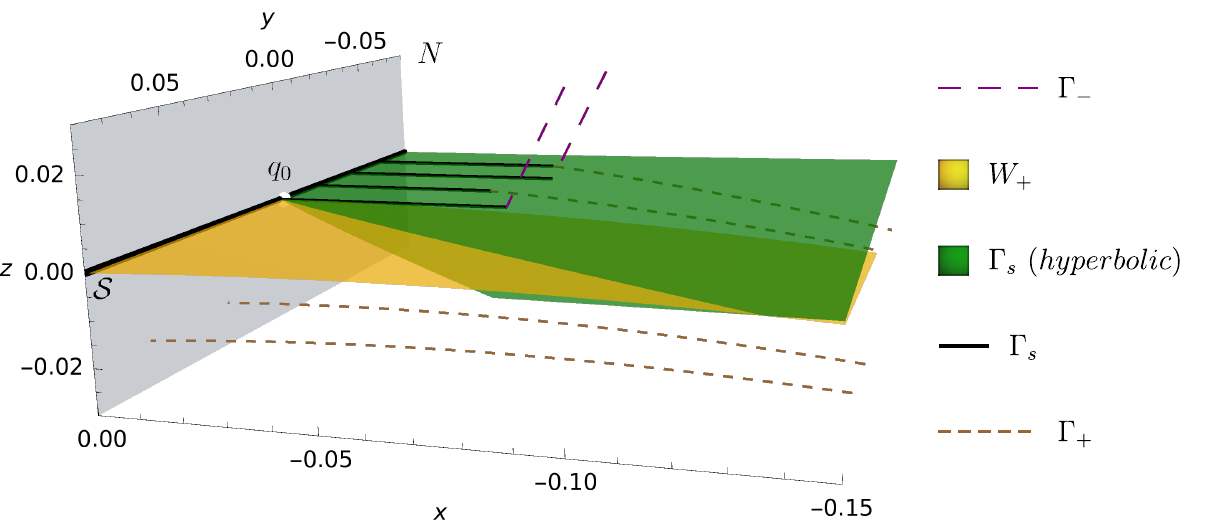}
\vspace*{0.5cm}

  \centering
  \begin{tabular}{ccc}
    \def\svgwidth{0.16\textwidth}
%% Accompanies image file 'saturate-a.pdf' (pdf, eps, ps)
\begingroup%
  \makeatletter%
  \providecommand\color[2][]{%
    \errmessage{(Inkscape) Color is used for the text in Inkscape, but the package 'color.sty' is not loaded}%
    \renewcommand\color[2][]{}%
  }%
  \providecommand\transparent[1]{%
    \errmessage{(Inkscape) Transparency is used (non-zero) for the text in Inkscape, but the package 'transparent.sty' is not loaded}%
    \renewcommand\transparent[1]{}%
  }%
  \providecommand\rotatebox[2]{#2}%
  \newcommand*\fsize{\dimexpr\f@size pt\relax}%
  \newcommand*\lineheight[1]{\fontsize{\fsize}{#1\fsize}\selectfont}%
  \ifx\svgwidth\undefined%
    \setlength{\unitlength}{152.45344162bp}%
    \ifx\svgscale\undefined%
      \relax%
    \else%
      \setlength{\unitlength}{\unitlength * \real{\svgscale}}%
    \fi%
  \else%
    \setlength{\unitlength}{\svgwidth}%
  \fi%
  \global\let\svgwidth\undefined%
  \global\let\svgscale\undefined%
  \makeatother%
  \begin{picture}(1,1.25068879)%
    \lineheight{1}%
    \setlength\tabcolsep{0pt}%
    \put(0,0){\includegraphics[width=\unitlength,page=1]{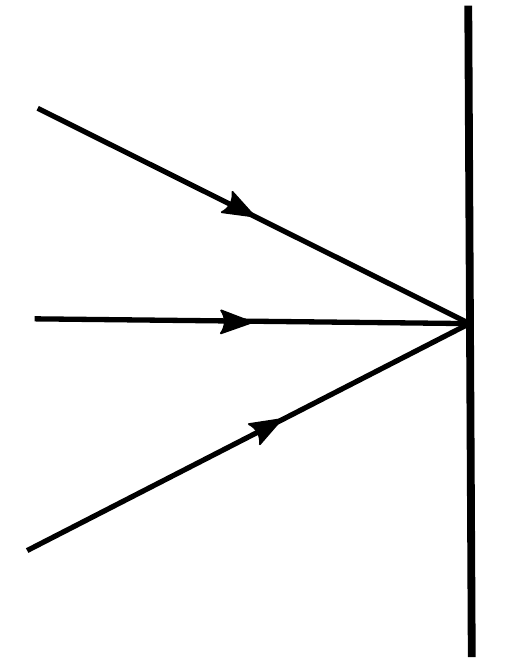}}%
    \put(0.10451818,1.07079054){\color[rgb]{0,0,0}\makebox(0,0)[lt]{\lineheight{1.25}\smash{\begin{tabular}[t]{l}{\tiny $\sigma_-$}\end{tabular}}}}%
    \put(0.93301918,1.00706186){\color[rgb]{0,0,0}\makebox(0,0)[lt]{\lineheight{1.25}\smash{\begin{tabular}[t]{l}$N$\end{tabular}}}}%
    \put(0.16240309,0.18966015){\color[rgb]{0,0,0}\makebox(0,0)[lt]{\lineheight{1.25}\smash{\begin{tabular}[t]{l}{\tiny $\sigma_+$}\end{tabular}}}}%
    \put(0.06252643,0.53721773){\color[rgb]{0,0,0}\makebox(0,0)[lt]{\lineheight{1.25}\smash{\begin{tabular}[t]{l}{\tiny $\sigma_s$}\end{tabular}}}}%
    \put(0,0){\includegraphics[width=\unitlength,page=2]{saturate-a.pdf}}%
  \end{picture}%
\endgroup%
 \hspace{0.2cm}
   &
    \def\svgwidth{0.29\textwidth}
   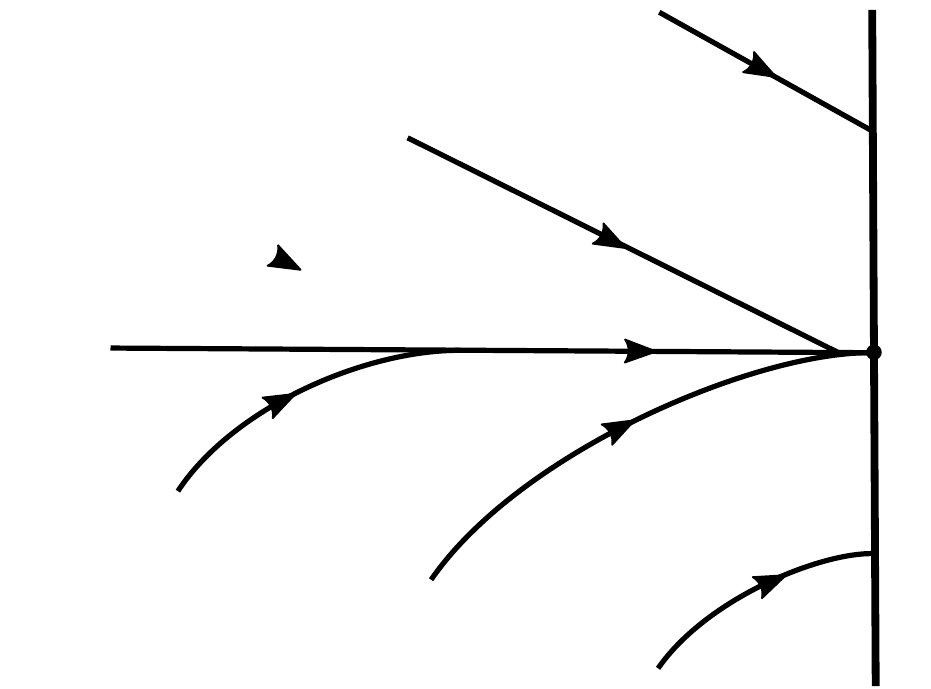 \hspace{0.2cm}
   &
    \def\svgwidth{0.33\textwidth}
   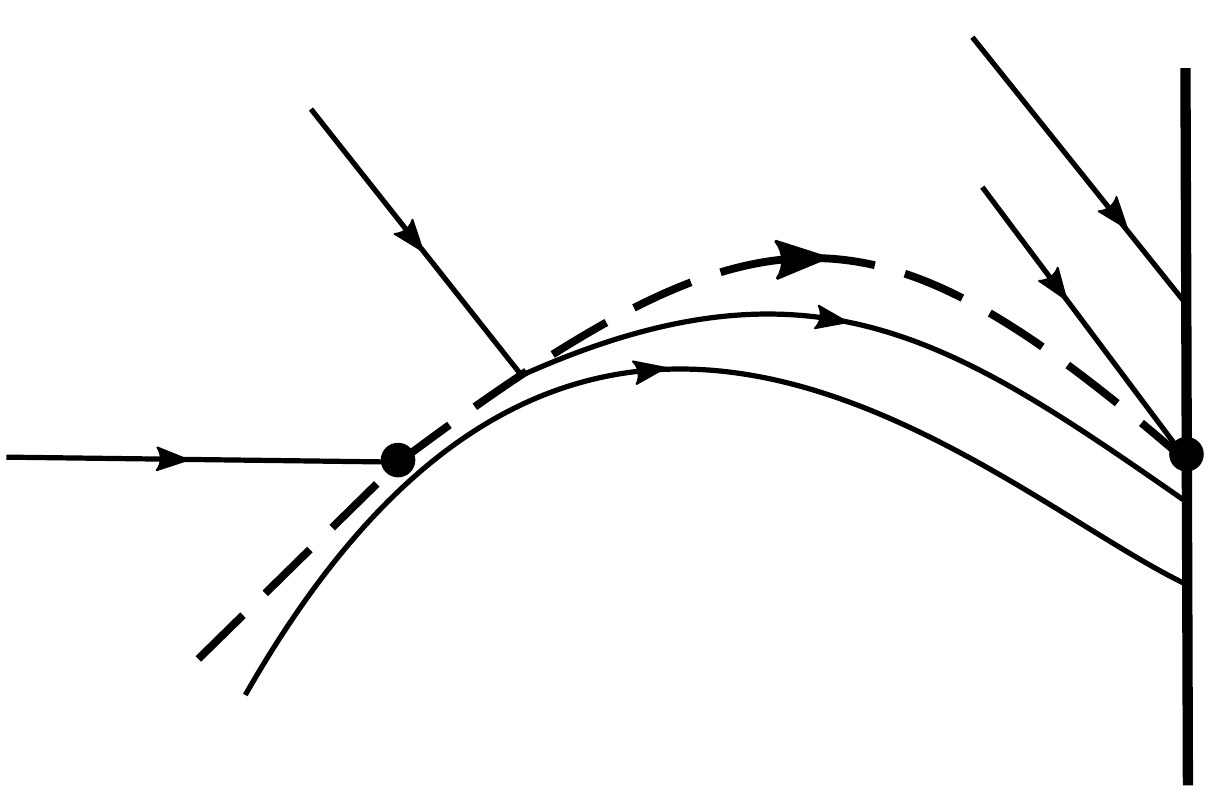
   \\
   $y<0$ & $y=0$ & $y>0$
     \end{tabular}
  \caption{Saturating case 1: $-1=a<0$, $1=u_{sx}>0$.
    Top figures are obtained using the \texttt{Mathematica}'s 
    program given in the Appendix \ref{prog:strata}
   and these
  situations are sketched in the bottom figures. 
\label{fig:saturation_case1}}
\end{figure}

\begin{figure}
  \centering
\includegraphics[scale=0.55]{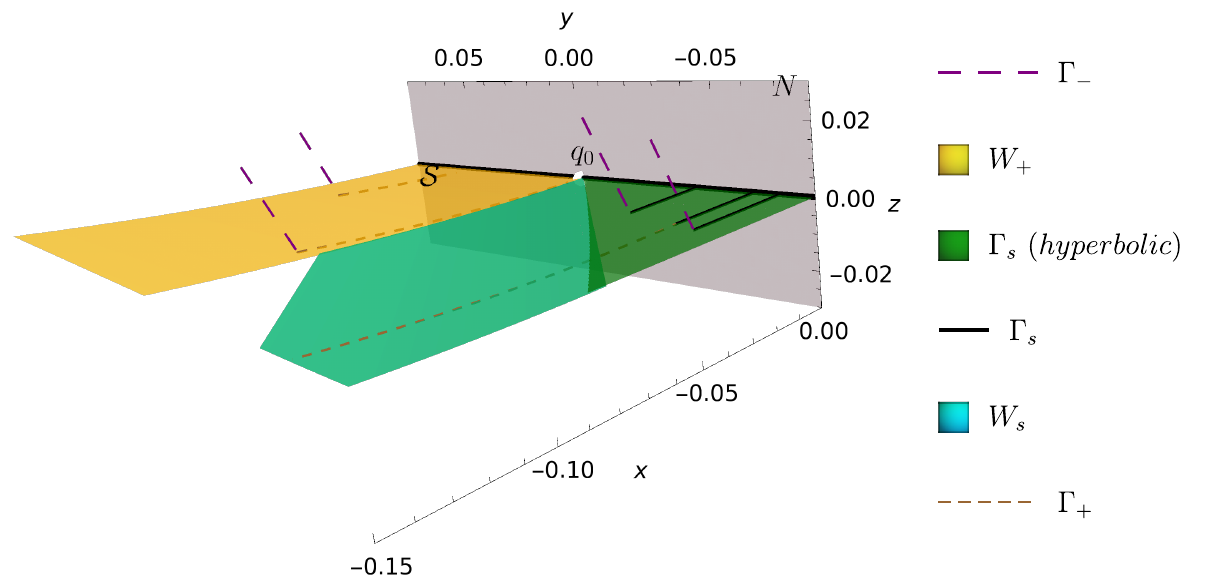} 
\vspace*{0.5cm}

    \begin{tabular}{ccc}
    \def\svgwidth{0.3\textwidth}
    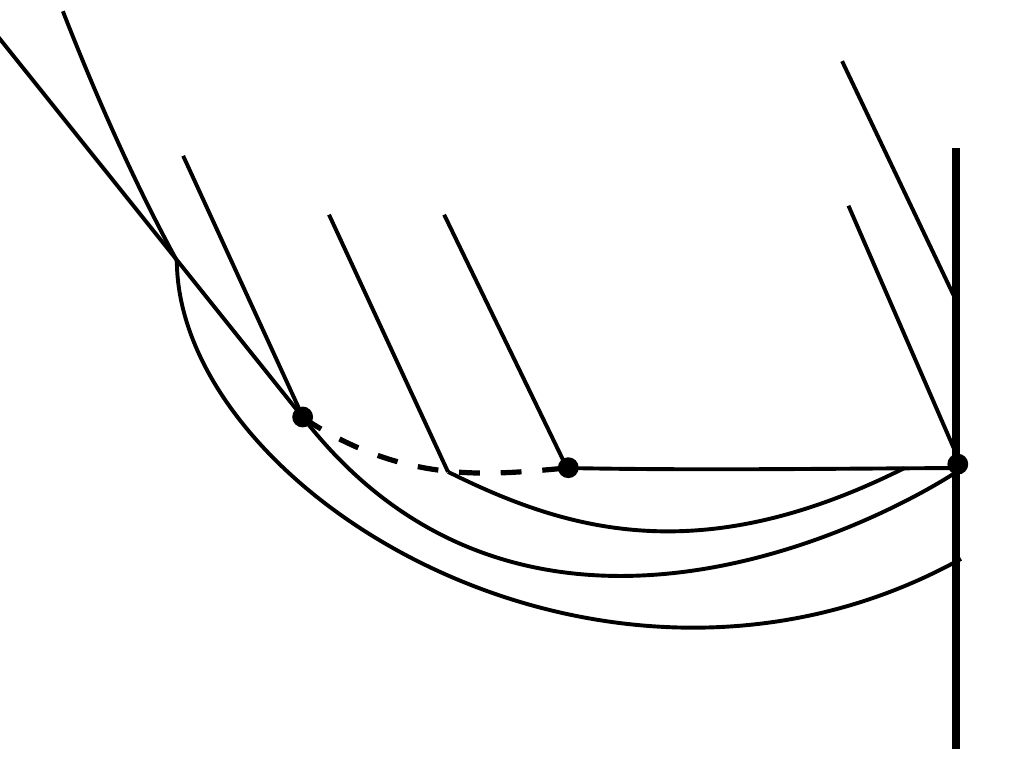 \hspace{0.1cm}
   &
    \def\svgwidth{0.26\textwidth}
%% Accompanies image file 'saturate-e.pdf' (pdf, eps, ps)
\begingroup%
  \makeatletter%
  \providecommand\color[2][]{%
    \errmessage{(Inkscape) Color is used for the text in Inkscape, but the package 'color.sty' is not loaded}%
    \renewcommand\color[2][]{}%
  }%
  \providecommand\transparent[1]{%
    \errmessage{(Inkscape) Transparency is used (non-zero) for the text in Inkscape, but the package 'transparent.sty' is not loaded}%
    \renewcommand\transparent[1]{}%
  }%
  \providecommand\rotatebox[2]{#2}%
  \newcommand*\fsize{\dimexpr\f@size pt\relax}%
  \newcommand*\lineheight[1]{\fontsize{\fsize}{#1\fsize}\selectfont}%
  \ifx\svgwidth\undefined%
    \setlength{\unitlength}{180.65055847bp}%
    \ifx\svgscale\undefined%
      \relax%
    \else%
      \setlength{\unitlength}{\unitlength * \real{\svgscale}}%
    \fi%
  \else%
    \setlength{\unitlength}{\svgwidth}%
  \fi%
  \global\let\svgwidth\undefined%
  \global\let\svgscale\undefined%
  \makeatother%
  \begin{picture}(1,0.82933707)%
    \lineheight{1}%
    \setlength\tabcolsep{0pt}%
    \put(0,0){\includegraphics[width=\unitlength,page=1]{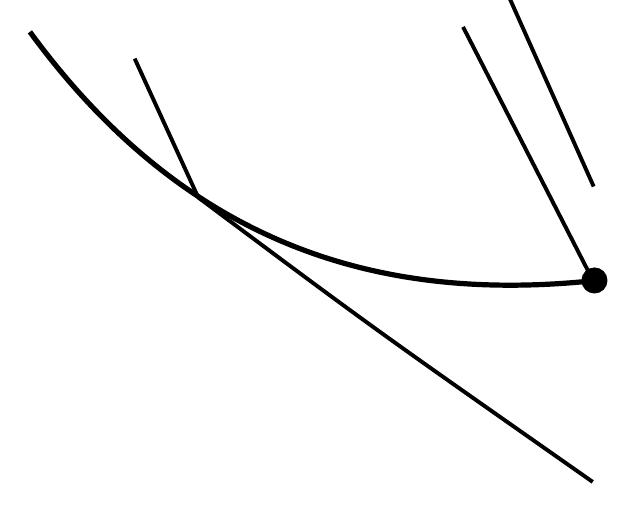}}%
    \put(0.17771508,0.76445706){\color[rgb]{0,0,0}\makebox(0,0)[lt]{\lineheight{1.25}\smash{\begin{tabular}[t]{l}{\tiny $\sigma_-$}\end{tabular}}}}%
    \put(0.67548044,0.81209455){\color[rgb]{0,0,0}\makebox(0,0)[lt]{\lineheight{1.25}\smash{\begin{tabular}[t]{l}{\tiny $\sigma_-$}\end{tabular}}}}%
    \put(0.61483165,0.71618487){\color[rgb]{0,0,0}\makebox(0,0)[lt]{\lineheight{1.25}\smash{\begin{tabular}[t]{l}{\tiny $\sigma_-$}\end{tabular}}}}%
    \put(0.97461205,0.72919614){\color[rgb]{0,0,0}\makebox(0,0)[lt]{\lineheight{1.25}\smash{\begin{tabular}[t]{l}$N$\end{tabular}}}}%
    \put(0.47346939,0.20419625){\color[rgb]{0,0,0}\makebox(0,0)[lt]{\lineheight{1.25}\smash{\begin{tabular}[t]{l}{\tiny $\sigma_+$}\end{tabular}}}}%
    \put(-0.00141045,0.60366721){\color[rgb]{0,0,0}\makebox(0,0)[lt]{\lineheight{1.25}\smash{\begin{tabular}[t]{l}{\tiny $W_+$}\end{tabular}}}}%
    \put(0,0){\includegraphics[width=\unitlength,page=2]{saturate-e.pdf}}%
  \end{picture}%
\endgroup%
\hspace{0.4cm}
   &
    \def\svgwidth{0.26\textwidth}
%% Accompanies image file 'saturate-f.pdf' (pdf, eps, ps)
\begingroup%
  \makeatletter%
  \providecommand\color[2][]{%
    \errmessage{(Inkscape) Color is used for the text in Inkscape, but the package 'color.sty' is not loaded}%
    \renewcommand\color[2][]{}%
  }%
  \providecommand\transparent[1]{%
    \errmessage{(Inkscape) Transparency is used (non-zero) for the text in Inkscape, but the package 'transparent.sty' is not loaded}%
    \renewcommand\transparent[1]{}%
  }%
  \providecommand\rotatebox[2]{#2}%
  \newcommand*\fsize{\dimexpr\f@size pt\relax}%
  \newcommand*\lineheight[1]{\fontsize{\fsize}{#1\fsize}\selectfont}%
  \ifx\svgwidth\undefined%
    \setlength{\unitlength}{194.58161926bp}%
    \ifx\svgscale\undefined%
      \relax%
    \else%
      \setlength{\unitlength}{\unitlength * \real{\svgscale}}%
    \fi%
  \else%
    \setlength{\unitlength}{\svgwidth}%
  \fi%
  \global\let\svgwidth\undefined%
  \global\let\svgscale\undefined%
  \makeatother%
  \begin{picture}(1,0.87037038)%
    \lineheight{1}%
    \setlength\tabcolsep{0pt}%
    \put(0,0){\includegraphics[width=\unitlength,page=1]{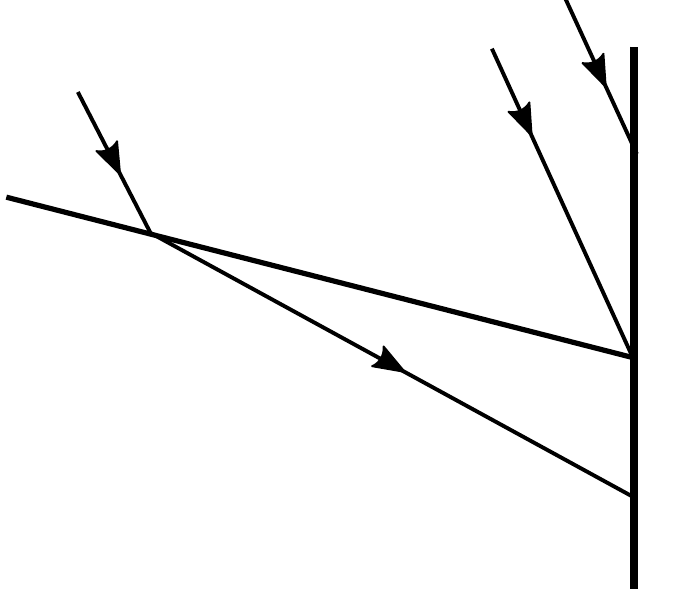}}%
    \put(0.06851852,0.74629631){\color[rgb]{0,0,0}\makebox(0,0)[lt]{\lineheight{1.25}\smash{\begin{tabular}[t]{l}{\tiny $\sigma_-$}\end{tabular}}}}%
    \put(0.68989459,0.85219267){\color[rgb]{0,0,0}\makebox(0,0)[lt]{\lineheight{1.25}\smash{\begin{tabular}[t]{l}{\tiny $\sigma_-$}\end{tabular}}}}%
    \put(0.67693154,0.57997046){\color[rgb]{0,0,0}\makebox(0,0)[lt]{\lineheight{1.25}\smash{\begin{tabular}[t]{l}{\tiny $\sigma_-$}\end{tabular}}}}%
    \put(0.37773715,0.31658067){\color[rgb]{0,0,0}\makebox(0,0)[lt]{\lineheight{1.25}\smash{\begin{tabular}[t]{l}{\tiny $\sigma_+$}\end{tabular}}}}%
    \put(-0.09899444,0.50404445){\color[rgb]{0,0,0}\makebox(0,0)[lt]{\lineheight{1.25}\smash{\begin{tabular}[t]{l}{\tiny $W_+$}\end{tabular}}}}%
    \put(0,0){\includegraphics[width=\unitlength,page=2]{saturate-f.pdf}}%
    \put(0.96768874,0.72674862){\color[rgb]{0,0,0}\makebox(0,0)[lt]{\lineheight{1.25}\smash{\begin{tabular}[t]{l}$N$\end{tabular}}}}%
  \end{picture}%
\endgroup%

     \end{tabular}
  \caption{Saturating case 2: $-1=a<0$, $1=u_{sx}<0$.
    Top figure is obtained using the \texttt{Mathematica}'s 
    program given in the Appendix \ref{prog:strata}
   and the corresponding situations 
   are sketched in the bottom figures. 
\label{fig:saturation_case2}}
\end{figure}

{\bf The case $\bm{a<0}$ and $\bm{u_s(0)=1}$.}
There are two generic cases described by Fig.\ref{fig:saturation_case1}-\ref{fig:saturation_case2}.
The two cases are discriminated by the existence or not of the singular arc 
in the transition. In the second case, the switching locus $\Sigma$ has two strata :
$\Sigma = W_+ \cup W_s$, where $W_+$ corresponds to optimal policies $\sigma_-\sigma_+$ and $W_s$
to policies $\sigma_-\sigma_+\sigma_s$, $\Sigma$ is not $C^1$. 
The stratum $W_s\cup \Gamma_s$ is $C^1$ but not $C^2$.


\subsubsection{Beyond the strict generalized Legendre--Clebsch condition}

In our previous work in the 90's we restricted our study to the case, where the 
strict generalized
Legendre--Clebsch condition is satisfied, see \cite{bonnard1997}. 
Now our program is to extend this analysis when this condition is not satisfied.

\paragraph*{The semi--bridge phenomenon as transition between two saturations}

A bridge in optimal control was introduced in \cite{cots2018} as a policy of the form
Bang-Singular-Bang-Singular, 
where the second bang is related as a saturation due to the violation 
of the strict Legendre-Clebsch condition.
In our context a semi--bridge is interpreted as a path to a bridge due to saturation
between optimal singular arcs.

{\bf A tutorial model for the semi--bridge.}
Consider the case 
\begin{equation}
\begin{aligned}
&\dot x =  1 + a\, y  -3c\,yz +c\, z^3 
\\
&\dot y =  z 
\\
&\dot z =  u, \qquad |u|\le 1,  
\end{aligned}
\label{eq:normal-form}
\end{equation}
where $a,c\neq 0$ are parameters.

{\it Lie brackets computations.}
We have
\[\begin{array}{ll}
F = (1+ay-3cyz+cz^3)\frac{\partial}{\partial x} + z \frac{\partial}{\partial y}, 
&
G = \frac{\partial}{\partial z},
\\
 \left[G,F\right] = 3 c (y-z^2) \frac{\partial}{\partial x} -  \frac{\partial}{\partial y},
&
[[G,F],G] = -6 c z \frac{\partial}{\partial x},
\\
\left[[G,F],F\right] = a\frac{\partial}{\partial x},
&
D(q) = \det(G,[G,F],[[G,F],G])=	-6 c z,
\\
D'(q)=\det(G,[G,F],[[G,F],F])= a,
&
D''(q) =  \det(G,[G,F],F)= a y-2 c z^3+1. 
\end{array}
\]

{\it Stratification of $N$.}\\
We have the following properties
\begin{itemize}
  \item $N:x=0$, $G=\frac{\partial}{\partial z}$
  \item ${\cal S}: n \cdot [G,F] = 0 \cap \{x=0\}$ is the parabola: $y=z^2$.
  \item $n \cdot \left[ \left[ G,F \right],G \right](0) =0$, 
    $n \cdot \left[ \left[ G,F \right],F \right](0) \neq 0$.
\end{itemize}
Observe that ${\cal E}: n \cdot F(q)=0$ is the defined by 
$y(a-3c z) + cz^3 +1 = 0$.

\begin{figure}
\centering
\includegraphics[scale=0.37]{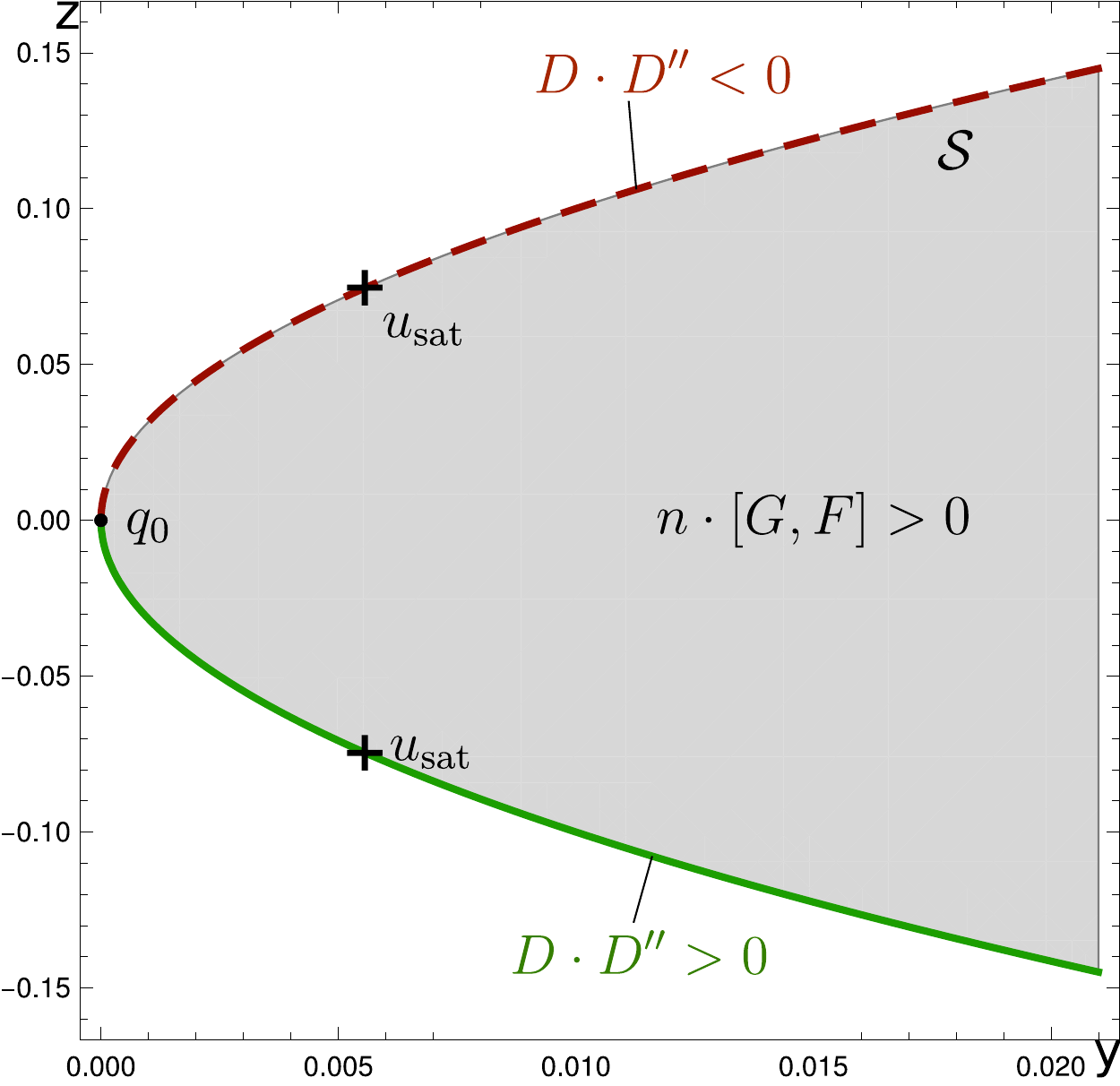}
\caption{Stratification of $N$ for the system \eqref{eq:normal-form} ($a=c=1$).
Green dotted line: elliptic, red line: hyperbolic, crosses: saturating values of 
the singular control.\label{fig:normal-form}}
\end{figure}

One will localize our study in a neighborhood $V$ of $0$ such that:
${\cal E} \cap V=\emptyset$. 
We represent such a situation on Fig.\ref{fig:normal-form}
One has $u_s = a/(6cz)$ and the points $u_{sat}$ are given by $|u_s|=1$.

 We have two saturating
points, one in the hyperbolic domain and one in the elliptic domain.
Near the fold, the two saturation phenomenons are glued together using the curvature
of ${\cal S}$, the normal to $N$ being given by $n=(1,0,0)$, while $q_0$
follows ${\cal S}$. Note that since the control is blowing up 
at the fold, we encounter the case $a>0, u_s(0)=3$ 
(see Fig.\ref{fig:par-ell})
and the bifurcation phenomenon of $C$, which splits into the case $\sigma_+$, $\sigma_-$ and the 
case $\sigma_+\sigma_-$ and $\sigma_-$ intersecting minimizers.
The two strata will be denoted $C_{12}$ and $C_1$. The model detects the two strata.

\begin{figure}[H]
    \centering
    \def\svgwidth{0.7\textwidth}
    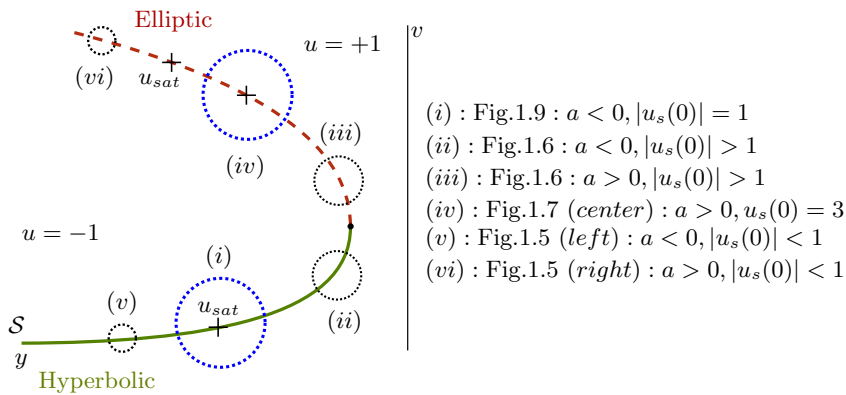
  \caption{Stratification of the target $x=0$ for a semi--bridge model and the
  associated local syntheses.}
  \label{fig:zoom}
\end{figure}

The stratification of the terminal manifold near the semi--bridge and the local time minimal synthesis
are described in Fig.\ref{fig:zoom}.

{\bf Local syntheses by symbolic computations.}

We use the \texttt{Mathematica}'s program in Appendix \ref{prog:strata} 
to derive the local syntheses for the semi--bridge model.

We compute the time minimal synthesis in a neighborhood $U$ of $q_0\in N$ where 
$N:x=0$ is the target.
Since the problem is flat, $q_0$ is either an ordinary switching point, a fold point, or
a codimension two case.

Denote $n$ the unit normal of $N$ at $q_0$ such that $n$ belongs to the half-space
containing the set $\{X+u\, Y,\ |u|\le 1\}$, i.e. $n=(1,0,0)$.

Assume $X(q_0)$ and $[Y,X](q_0)$ are not tangent to $N$. 
Then $q_0$ is an ordinary point and 
the optimal control is given by $u(q_0)=-\text{sign} \ p\cdot [Y,X](q_0)$. 

We take $q_0 = (0,s_0^2,s_0)$  and we look at the behavior of 
BC-extremals reaching $N$ near $q_0$.
The switching surface $W$, the splitting locus $C$ and the trajectories $\sigma_\pm$ 
are computed via symbolic computations, expanding the jets of $F$ and $G$ up to some order.

\begin{figure}[H]
\centering
\includegraphics[scale=0.6]{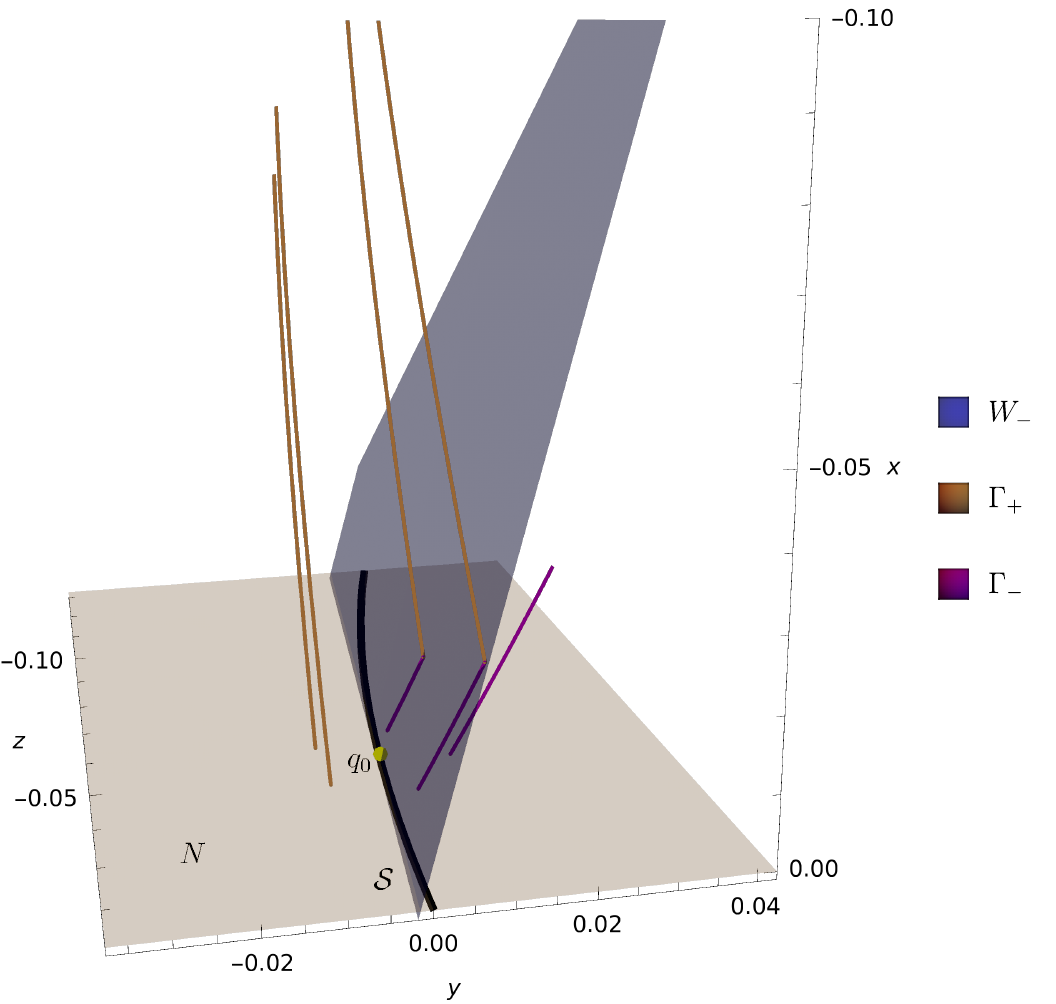}
  \caption{Time minimal synthesis near an hyperbolic fold point for $q_0=(0,(-0.05)^2,-0.05)$. Only the BC-extremals
  with optimal control $-1$ switch for the tutorial model \eqref{eq:normal-form} with $a=c=1$.
The tangent space of the switching surface at $q_0$ is represented.  
   \label{fig:nf-near-hyperbolic}  
}
\end{figure}

\begin{figure}
\centering
\includegraphics[height=5.2cm]{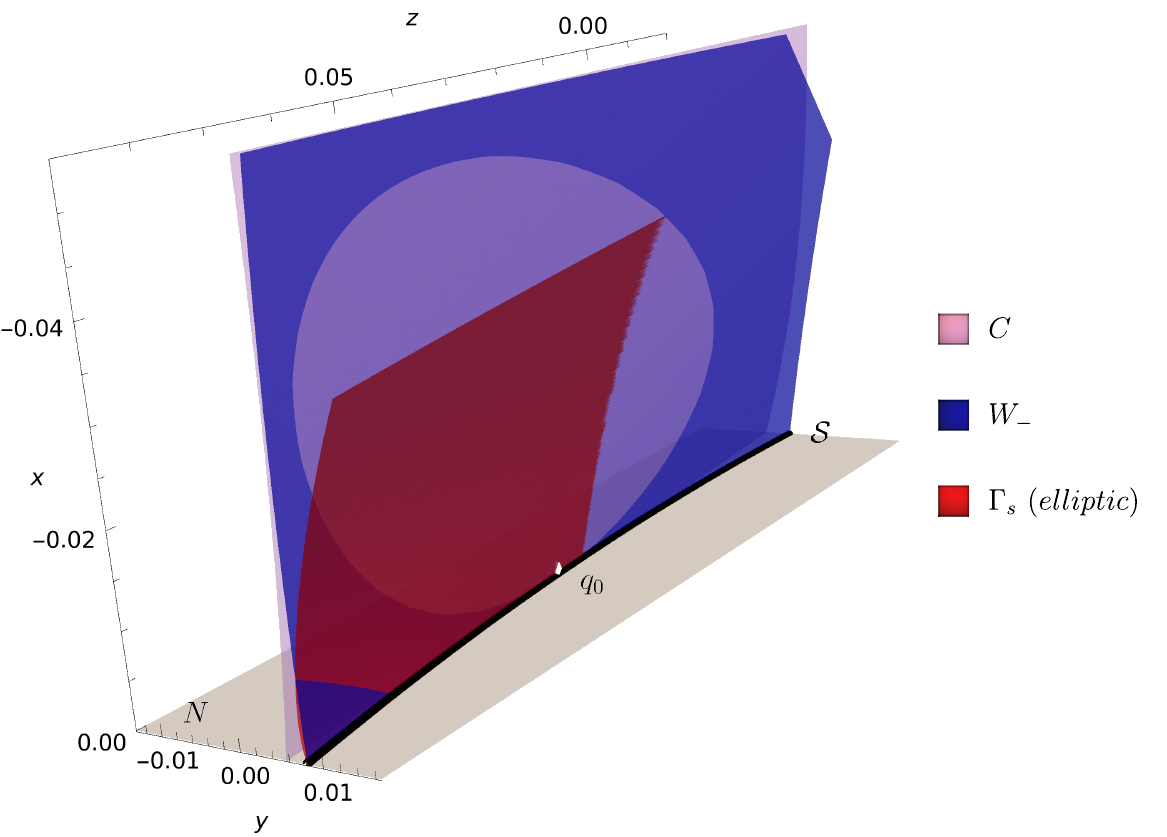}
  \caption{Time minimal synthesis near an elliptic fold point for $q_0=(0,0.04^2,0.04)$ 
  for the tutorial model \eqref{eq:normal-form} with $a=1,\,c=5$.
   \label{fig:nf-near-elliptic}  
}
\end{figure}

\begin{figure}
\centering
\includegraphics[height=5.2cm]{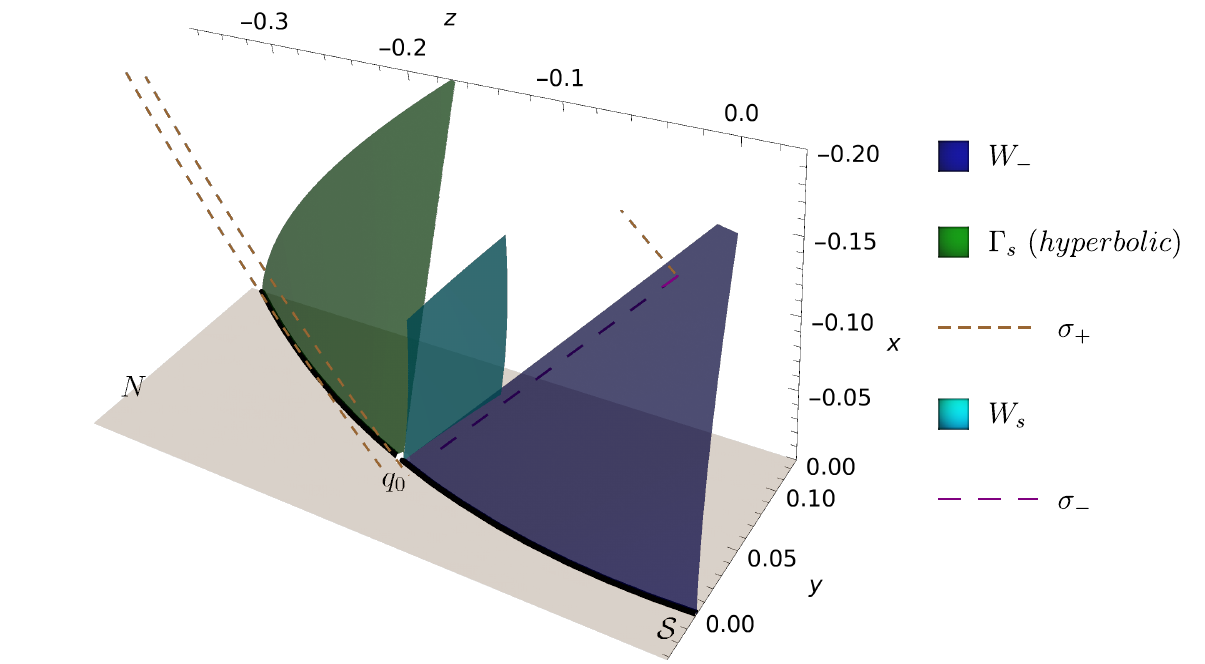}
\caption{Time minimal synthesis near a saturating hyperbolic fold point  $q_0=(0,z_{sat}^2,-z_{sat})$
  for the tutorial model \eqref{eq:normal-form} with $a=c=1$. 
  Trajectories $\sigma_-$ switch, while $\sigma_+$ don't. The surface $W_s$ is the first 
  switching point of the sequence $\sigma_+\sigma_-\sigma_s$.
   \label{fig:nf-hyp-sat}  
}
\end{figure}

\begin{itemize}
\item
\underline{Near a parabolic point on an hyperbolic singular arc.}\\
BC-extremals $\sigma_-$ switch while BC-extremals $\sigma_+$ don't.
These computations are represented in Fig. \ref{fig:nf-near-hyperbolic}
for $q_0=(0,(-0.05)^2,-0.05)$. 
\item
\underline{Near a parabolic point on an elliptic singular arc.}\\
BC-extremals $\sigma_\pm$ reflect on the switching surface $W_\pm$, hence a 
cut locus appears between the switching surfaces.
We represent in Fig. \ref{fig:nf-near-elliptic} with $q_0=(0,0.04^2,0.04),\, 
0.04>z_{sat}=a/(6c)$ this cut locus 
and the switching surface $W_-$, together with the
elliptic singular surface $\Gamma_s$.  
\item
\underline{Near a saturated point on an hyperbolic singular arc.}\\
We have $q_0=(0,z_{sat}^2,-z_{sat})$ where $z_{sat}=\frac{a}{6c}$.
The synthesis is represented in Fig.\ref{fig:nf-hyp-sat} and it corresponds
to the synthesis described earlier by Fig.\ref{fig:saturation_case1}.
\end{itemize}

{\it These computations confirm the expected results and validate the symbolic program to investigate
the McKeithan model.}

%\vspace*{2cm}
\begin{remark}
We can investigate other singularities by the same method, for instance:
\begin{itemize}
  \item Cusp singularity $y=z^3$.
  \item Singularity $y^2=z^2$.
  \item Case $D=D'=0$ is $q=\{0\}$.
\end{itemize}
\end{remark}

{\bf Global analysis.}

Take $u=\varepsilon\in \{-1,1\}$. The dynamics \eqref{eq:normal-form} is integrable and we have
\[
  p_3(t) =   \frac{1}{2} t \left(-2 c t^2+ (a -6 c \varepsilon s_0) t-6 c s_0^2+6
  c w_0\right).
\]

Denoting by $\Gamma_\varepsilon(q_0),\ \varepsilon\in\{-1,1\}$ the surface composed by the 
trajectories $\sigma_\varepsilon$ passing through $q_0$, then $\Gamma_\varepsilon(0,w_0,s_0)$ 
is parameterized by
\begin{equation*}
  \begin{aligned}
    x = & t \left( (a -3cs_0) + c s_0^3+1\right) + 
  \frac{1}{2} t^2 \left(s_0 (a-3 c s_0)+3
  c\varepsilon  (s_0^2- w_0)\right) + o(t^3)
    \\
  y = &w_0 + s_0 t+\frac{\varepsilon t^2}{2}+ o(t^3) 
  \\
  z = &s_0+ \varepsilon t + o(t^3). 
  \end{aligned}
  \label{eq:nf-gamma}
\end{equation*}

\begin{figure}[htpb]
  \centering
  \includegraphics[scale=0.62]{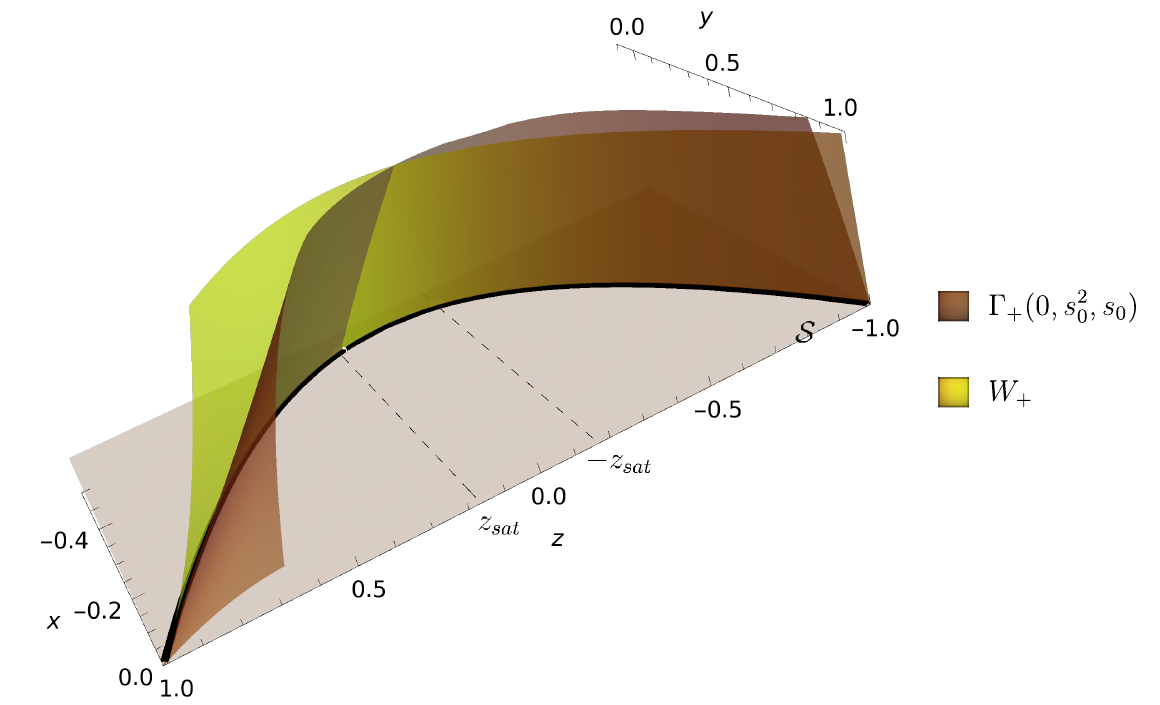}
  \caption{The trajectories $\sigma_+$ starting from $q_0=(0,w_0,s_0)$ such that $w_0<s_0^2$ 
  switch on $W_+$ for $z\in ]z_{sat},1[$ ($a=c=1$).}
  \label{fig:wp-cap-gp}
\end{figure}

\begin{figure}[htpb]
  \centering
  \includegraphics[scale=0.62]{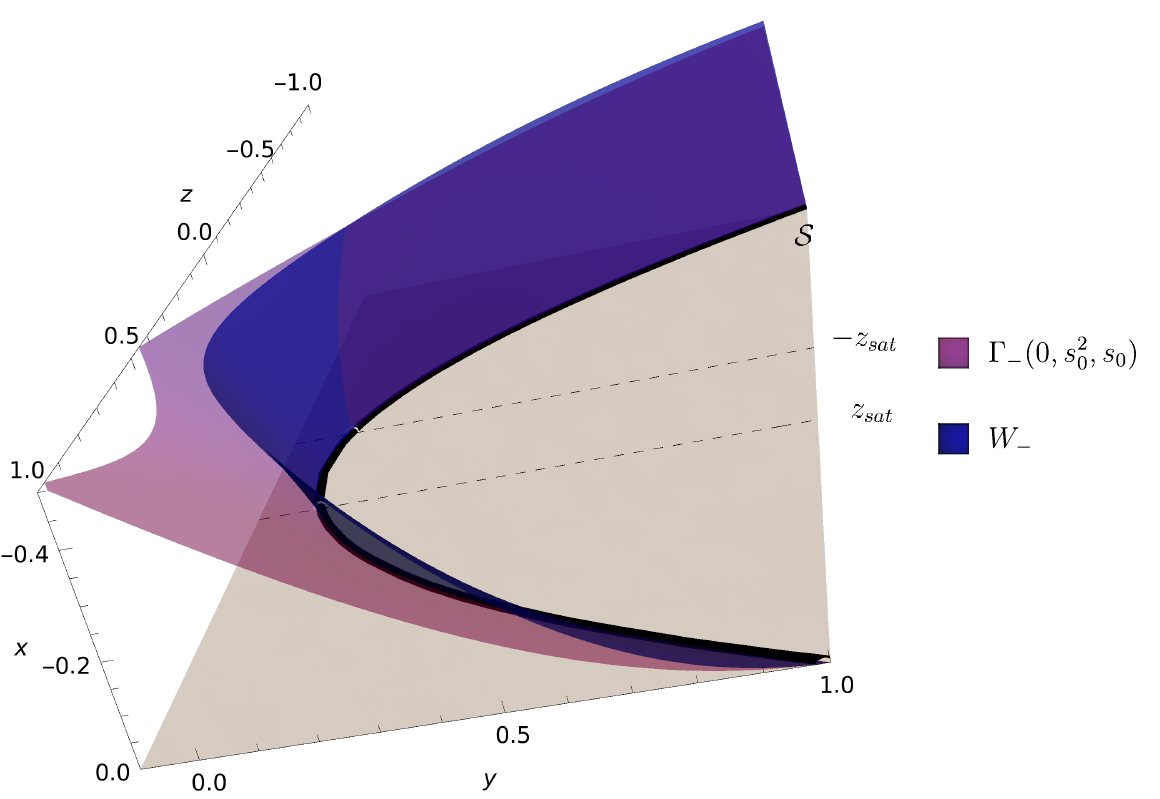}
  \caption{The trajectories $\sigma_-$ starting from  $q_0=(0,w_0,s_0)$ such that $w_0>s_0^2$ 
  switch on $W_-$ for $z\in ]-z_{sat},1[$ ($a=c=1$).}
  \label{fig:wm-cap-gm}
\end{figure}

Solving $p_3=0$ with respect to $w_0$, the switching surface $K_\varepsilon$ 
for BC-extremals of controle $u=\varepsilon$ on $N$
is parameterized by:
\begin{equation*}
  \begin{aligned}
  x = & \frac{1}{6} t^2 \left(-\frac{a^2}{c}+6 a \varepsilon s_0+6 a
    s_0-18 c \varepsilon s_0^2-9 c s_0^2\right)+\frac{1}{6} t
    \left(6 a s_0^2-12 c s_0^3+6\right),
    \\
  y = & t \left((\varepsilon+1) s_0-\frac{a}{6 c}\right)
  +\frac{1}{6} (3 \varepsilon+2) t^2+s_0^2 
  \\
  z = &s_0+ \varepsilon t + o(t^3). 
  \end{aligned}
  \label{eq:nf-w}
\end{equation*}
The switching condition of a trajectory on $W_\varepsilon$
is given by the sign of
\[
\det\left(
\frac{\partial K_\varepsilon}{\partial t}_{\mid t=0},
\frac{\partial K_\varepsilon}{\partial t}_{\mid t=0},
\frac{\partial \Gamma_\varepsilon(0,s_0^2,s_0)}{\partial t}_{\mid t=0}
\right),
\]
and this determinant is equal to 
   $\left(a s_0^2-2 c s_0^3+1\right) (6 c \varepsilon
     s_0-a)/(6 c)$.
In the case $a=c=1$, we get that $\sigma_+$ (resp. $\sigma_-$) switches on $W_+$ (resp. $W_-$) 
for $z\in ]z_{sat},1[$ (resp. $z\in ]-z_{sat},1[$) as illustrated in Fig.\ref{fig:wp-cap-gp}
(resp. Fig.\ref{fig:wm-cap-gm}).

\begin{figure}[htpb]
  \centering
  \begin{tabular}{cc}
    \includegraphics[scale=0.55]{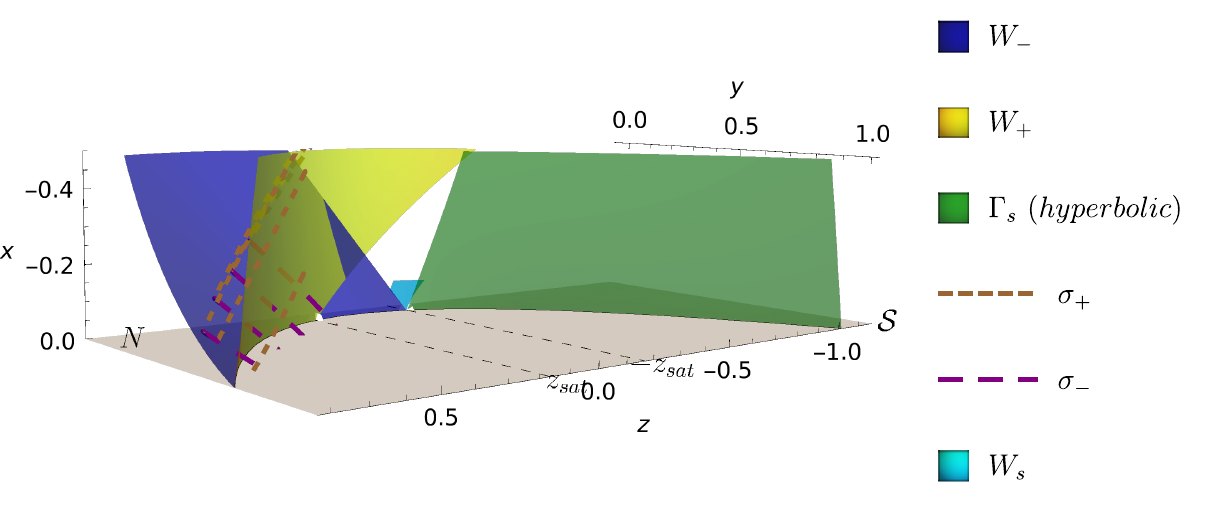}
& 
  \includegraphics[scale=0.55]{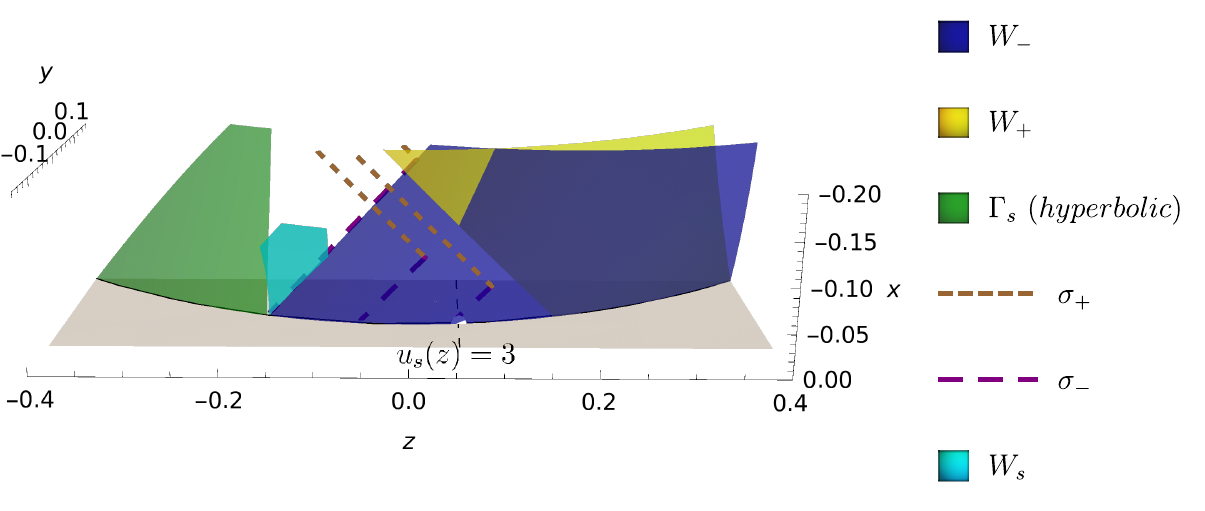}
\\
Cut locus for $z>z_{sat}$. &  $u_s(z)=3$.
  \end{tabular}
  \caption{ {\it (left)} $\sigma_-$ reflect on $W_-$ for $z>z_{sat}$. 
  {\it (right)} At $z=\frac{a}{18c}<z_{sat}$ such that $u_s(z)=3$, $\sigma_-$ cross $W_-$ and 
  the cut locus disappears (see Fig.\ref{fig:par-ell}). For $-z_{sat}<z<z_{sat}$, the optimal
  policy is $\sigma_+\sigma_-$. The case $z=-z_{sat}$ was given in Fig.\ref{fig:nf-hyp-sat}
  and is also represented here with the surface $W_s$.
  }
  \label{fig:nf-global-a}
\end{figure}

\begin{figure}[htpb]
  \centering
  \includegraphics[scale=0.62]{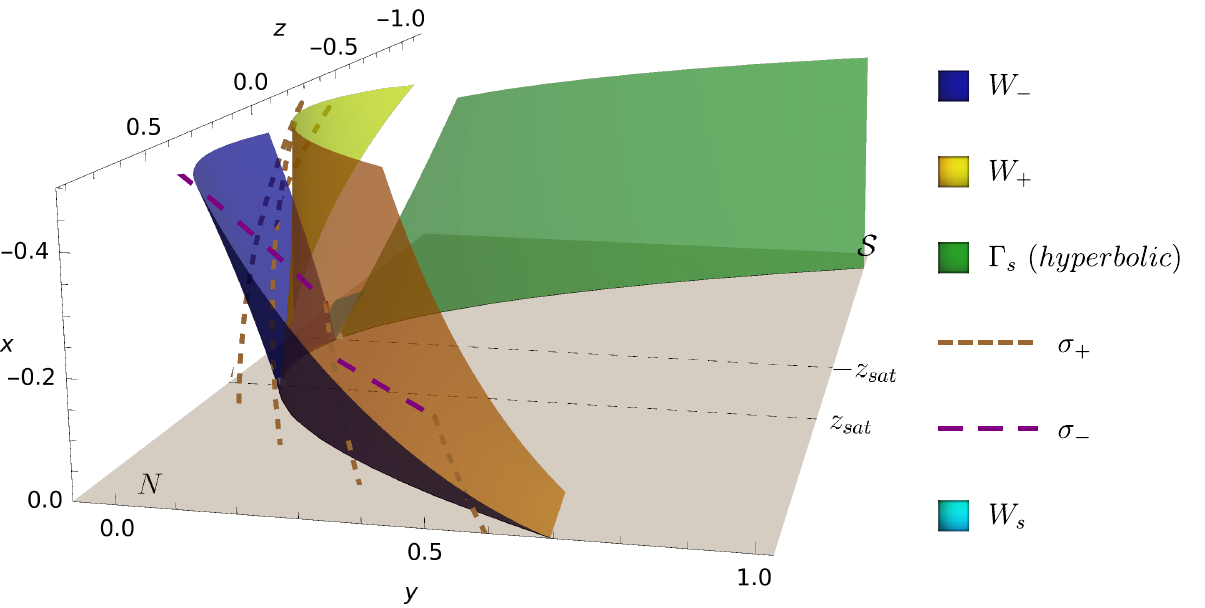}
  \caption{For $z>z_{sat}$ $\sigma_+$ reflect on $W_+$, hence we have a cut locus.
  For $z<-z_{sat}$, the optimal policy is Bang-Singular.
  }
  \label{fig:nf-global-b}
\end{figure}

We deduce the time minimal synthesis from Fig.\ref{fig:nf-global-a}-\ref{fig:nf-global-b}
for all point $q_0$ along ${\cal S}$.
More precisely, if $z<-z_{sat}$ the optimal policy is Bang-Singular-Bang;  
if $z=-z_{sat}$, the synthesis is given by Fig.\ref{fig:nf-hyp-sat}.
If $-z_{sat}<z<z_{sat}$, we have first a $\sigma_+\sigma_-$'s policy, then 
we encounter the case $u_s(z)=3$ (see Fig.\ref{fig:par-ell}), then we have a cut 
for $z>z_{sat}$.

\paragraph*{Construction of a semi--normal form and the need of integrability assumption.}

\noindent
{\it Notations and assumption:} We normalize in a small neighborhood of $N$.
We assume that $F$ is transverse to $N$.
Let $q=(x,y,z)$, $q_0=(0,0,0)$, $G=\frac{\partial}{\partial z}$ and $N\coloneqq \{x=0\}$.
We denote by $\sigma_0$ the reference trajectory of $F$ identified to 
$t\rightarrow (t,0,0)$ so that $F_{\mid \sigma_0} = \frac{\partial}{\partial x}$.

Assuming $x$ small enough one can write the system as:
\begin{align*}
  &\dot x = 1 + f(y,z) + R_1
  \\
  &\dot y = g(y,z) + R_2
  \\
  &\dot z = h(y,z) + R_3,
\end{align*}
where $R_i(x,y,z)=o(|x|)$ and the corresponding
model is 
\begin{align*}
  &\dot x = 1 + f(y,z)
  \\
  &\dot y = g(y,z)
  \\
  &\dot z = h(y,z) + u.
\end{align*}

\noindent
{\it Singular flow on the model:} Furthermore, we assume
$\frac{\partial g}{\partial z} \neq 0$. Since the singular flow is feedback invariant,
setting $z=g(y,z)$, $z\rightarrow z$ and using a proper feedback the system takes the form
\begin{align*}
  &\dot x = 1 + f(y,z) 
  \\
  &\dot y = z
  \\
  &\dot z = u.
\end{align*}
Computing, we get:
\begin{equation}
  \begin{array}{ll}
    \left[ G,F \right] = -f_z \frac{\partial}{\partial x} + \frac{\partial }{\partial y}
    & \left[ \left[ G,F \right],G \right] = -f_{z^2} \frac{\partial}{\partial x}, 
    \\
    \left[ \left[ G,F \right],F \right] = -f_{yz} \frac{\partial}{\partial x}
    +f_y \frac{\partial}{\partial y}. 
  \end{array}
\end{equation}
Hence
\[
  D = -f_{z^2}, \qquad D' = -f_y f_z  + f_{yz}
\]
and 
\[
  u_s = -\frac{D'}{D} = \frac{f_yf_z - f_{yz}}{f_{z^2}}.
\]
Therefore, we have
\[
  \frac{\dd y}{\dd z} = \frac{z}{u_s}
\]
and assuming the {\it separability condition}
\[
  \frac{\dd y}{\dd z} = h_1(y) h_2(z)
\]
one can integrate the singular trajectories such that $q(0) \in {\cal S}$ as follows, 
parameterizing by $z$ instead of $t$ using the relation
\[
  \int_{y(0)}^y h_1(y)\,\dd y = \int_{z(0)}^z h_2(z)\, \dd z.
\]

Assuming that $y(0)$ can be expressed as $y(0)=p_1(z(0))$ and moreover that
the left--hand side is invertible, one gets 
$y = p_2(z,z(0))$. This defines the {\it singular leaf} ${\cal F}$ through 
$q(0) \in {\cal S}$ using the equation
\[
  \frac{\dd x}{\dd z} = \frac{1+f(y,z)}{z}.
\]
Using $x(0)=0$, one gets:
\[
  x(z,z(0)) = \int_{z(0)}^z \left( 1 +f(p_2(z,z(0)),z) \right)\, \dd z.
\]
\begin{remark}
  This requires an integrability condition to obtain the leaf ${\cal F}$. 
  If it is not satisfied, integration is obtained through a proper {\it numeric
  integration}.
\end{remark}

\paragraph*{Application.} 

\begin{figure} 
\centering
\includegraphics[scale=0.64]{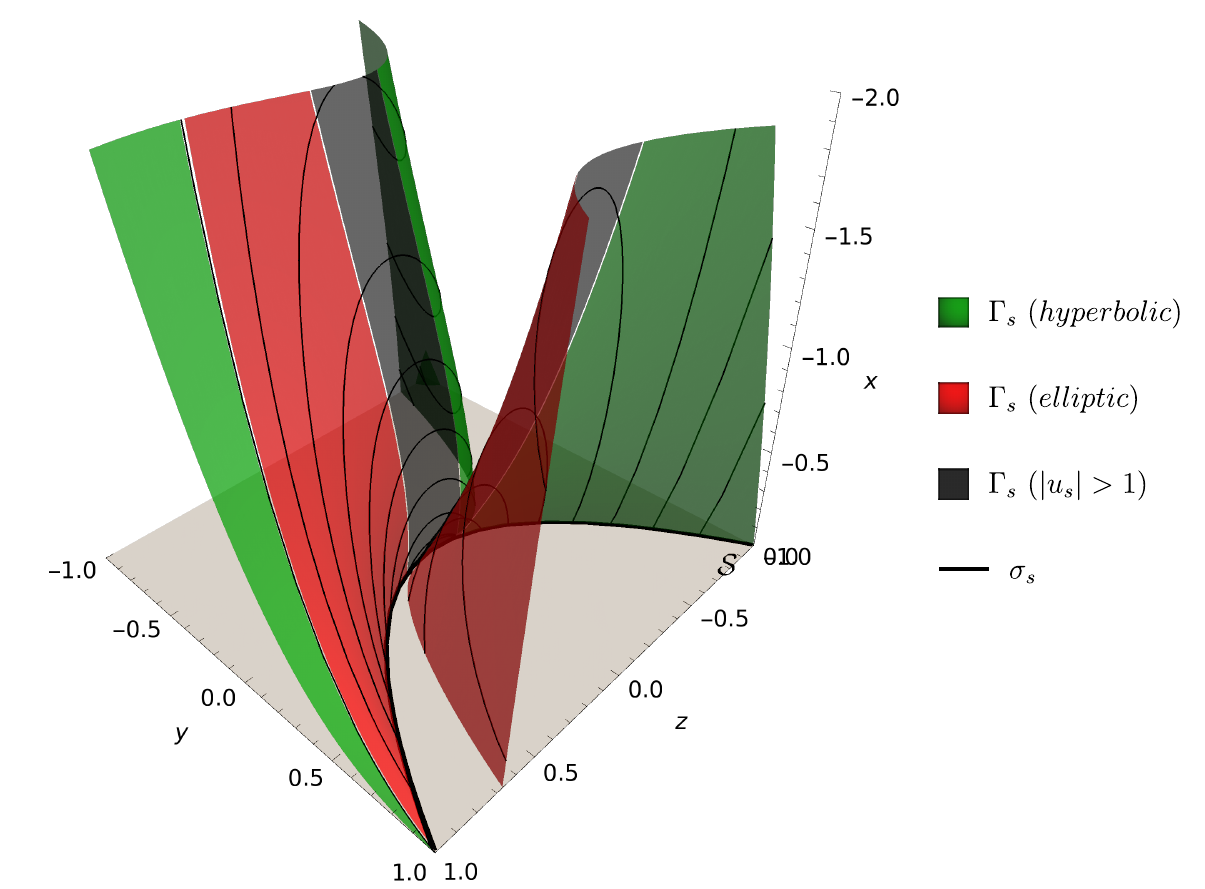}
\caption{Stratification of the singular set $\Gamma_s$ for the system \eqref{eq:normal-form}
with parameters $a=c=1$.\label{fig:sing-leaf}}
\end{figure}

{\bf Singular set.}
Back to the tutorial model \eqref{eq:normal-form}, we fix $a=c=1$,
we are in the separable case and using our previous algorithm one gets
\begin{proposition}
  \begin{enumerate}
    \item Integrating the singular flow, we get:
\[
\begin{array}{ll}
x(t,z(0))&= t + \frac{2\varepsilon (2a+3ck) \,\left(k t+z(0)^2\right)^{5/2}
-10 c \left(k t+z(0)^2\right)^3 t}{15 k^2}
\\ 
y(t,z(0))&= \frac{2 \varepsilon }{3k} \left(k t+z(0)^2\right)^{3/2}
\\
z(t,z(0))&= \varepsilon\, \left(k t+z(0)^2\right)^{1/2},
\end{array}
\]
where $\varepsilon=\pm 1$, $k=12c/a$.
\item 
  Parameterizing by $z$, the integral singular leaf ${\cal F}$ with $q(0)\in {\cal S}$
  is given
\[
  \begin{array}{l}
  y(z,z(0))=2c\, (z^3-z(0)^3)/a + z(0)^2  
  \\
  x(z,z(0)) = g(z)  - g(z(0)),
  \end{array}
\]
where 
$
\displaystyle
     g(z) = \frac{3 c z^2}{5a^2}\, 
     \Big(
       2 a (2 c+1)\, z^3+30 c z(0)^2 (a-2 c z(0))\, z 
   +30 c^2\, z^4  +5 a (a z(0)^2-2 c z(0)^3+1)\Big).
$
\end{enumerate}
\end{proposition}
We represent on Fig.\ref{fig:sing-leaf} the stratification of the singular leaf ${\cal F}$
with $a=c=1$ using hyperbolic and elliptic cases and the admissibility conditions
$|u_s|\le 1$.

\subsection{The 2d--McKeithan scheme}
Recall that the network is
$
\tikz[baseline={(a.base)},node distance=6mm] {
\node (a) {$T+M$};
\node[right=of a] (b) {$A$};
\node[right=of b] (c) {$B$};
\draw (a.base east) edge[left,->] node[above] {{$\scriptstyle k_1$}} (b.base west) (b.base east) edge[->] node[above] {{$\scriptstyle k_2$}} (c.base west);
\draw (a) edge[bend right=40,<-] node[below right] {{$\scriptstyle k_3$}} (b);
\draw (a.south) edge[bend right=60,<-] node[below right] {{$\scriptstyle k_4$}} (c);
}
$
and the system $\dot q=F(q)+u\, G(q)$, using the coordinates $q=(x,y,v)$,
$x=[A]$, $y=[B]$, $v=k_1$ and reduced to the stoichiometric class
$T+M+A=\delta_1$ and $T+M+B=\delta_2$ takes the form
\begin{align*}
&\dot x =    -\beta_2 x v^{\alpha_2}-\beta_3 x v^{\alpha_3}-\delta_3 v\, (x+y)+\delta_4 v+v\, (x+y)^2
\\
&\dot y = \beta_2 x v^{\alpha_2}-\beta_4 y v^{\alpha_4}
\end{align*}
with 
$0\le x\le \delta_1$, 
$0\le y \le \delta_2$, 
$\delta_3=\delta_1+\delta_2$,
$\delta_4=\delta_1\, \delta_2$, 
$k_2=\beta_2\, v^{\alpha_2}$, 
$k_3=\beta_3\, v^{\alpha_3}$,
$k_4=\beta_4\, v^{\alpha_4}$.

\subsubsection{Lie brackets computations}
We have:	
\begin{itemize}
\item 
$F = 
\big(-\beta_2 x v^{\alpha_2}-\beta_3 x v^{\alpha_3}-\delta_3 v (x+y)
+ \delta_4v+v\, (x+y)^2\big)\frac{\partial}{\partial x}
+
\big(\beta_2 x v^{\alpha_2}-\beta_4 y v^{\alpha_4}\big)\frac{\partial}{\partial y}
$,
\item 
$G = \frac{\partial}{\partial v}$,
\item 
$\left[ G,F \right]=
\big(x (\alpha_2 \beta_2 v^{\alpha_2-1}+\alpha_3 \beta_3
v^{\alpha_3-1}+\delta_3)+\delta_3 y-\delta_4-x^2-2 x
y-y^2\big)\frac{\partial}{\partial x}+
\big(\alpha_4 \beta_4 y v^{\alpha_4-1}-\alpha_2 \beta_2 x
v^{\alpha_2-1}\big) \frac{\partial}{\partial y}
$,
\item 
$
\left[ \left[ G,F \right],G \right]=
\big(x (\alpha_2 \beta_2(\alpha_2 -1)
v^{\alpha_2-2}+\alpha_3 \beta_3(\alpha_3-1)
v^{\alpha_3-2})\big)\frac{\partial}{\partial x}
+\big(y \alpha_4 \beta_4(\alpha_4 -1) v^{\alpha_4-2}
-x \alpha_2 \beta_2(\alpha_2-1) v^{\alpha_2-2}
\big)\frac{\partial}{\partial y}
$,
\item 
$
\left[ \left[ G,F \right]F \right]=
\big(-x v^{\alpha_2} \beta_2 \delta_3(\alpha_2-1)
-y v^{\alpha_2} \beta_2 \delta_3(\alpha_2 -1)
+v^{\alpha_2} (\alpha_2 \beta_2 \delta_4
+x y (2 \alpha_2 \beta_2-2 \beta_2)-\beta_2 \delta_4)
+x^2 (\alpha_2 \beta_2-\beta_2) v^{\alpha_2}
+y^2 (\alpha_2 \beta_2-\beta_2)v^{\alpha_2}
+y v^{\alpha_3} (\beta_3 \delta_3-\alpha_3 \beta_3
\delta_3)+v^{\alpha_3} (\alpha_3 \beta_3 \delta_4-\beta_3
\delta_4)+x^2 (\beta_3-\alpha_3 \beta_3) v^{\alpha_3}+y^2
(\alpha_3 \beta_3-\beta_3) v^{\alpha_3}+y v^{\alpha_4} (\alpha_4
\beta_4 \delta_3-\beta_4 \delta_3)+x y (2 \beta_4-2 \alpha_4
\beta_4) v^{\alpha_4}+y^2 (2 \beta_4-2 \alpha_4 \beta_4)
v^{\alpha_4}\big) \frac{\partial}{\partial x}
+\big(
x (v^{\alpha_2+\alpha_3-1} (\alpha_2 \beta_2
\beta_3-\alpha_3 \beta_2 \beta_3)+v^{\alpha_2+\alpha_4-1}
(\alpha_4 \beta_2 \beta_4-\alpha_2 \beta_2 \beta_4))+x
v^{\alpha_2} (\alpha_2 \beta_2 \delta_3-\beta_2 \delta_3)+y
v^{\alpha_2} (\alpha_2 \beta_2 \delta_3-\beta_2
\delta_3)+v^{\alpha_2} (-\alpha_2 \beta_2 \delta_4+x y (2
\beta_2-2 \alpha_2 \beta_2)+\beta_2 \delta_4)+x^2
(\beta_2-\alpha_2 \beta_2) v^{\alpha_2}+y^2 (\beta_2-\alpha_2
\beta_2) v^{\alpha_2}\big) \frac{\partial}{\partial y}
$.
\end{itemize}

\subsubsection{Singular arcs}

One has:
\begin{itemize}
\item
$D(q)=
((\alpha_4-1) \alpha_4 \beta_4 y v^{\alpha_4-3}
-(\alpha_2-1) \alpha_2 \beta_2 x v^{\alpha_2-3}
)
(\alpha_2 \beta_2 x v^{\alpha_2}+\alpha_3 \beta_3 x
v^{\alpha_3}+\delta_3 v (x+y)-\delta_4 v-v x^2-2 v x y-v y^2)
+x
(\alpha_2^2 \beta_2 v^{\alpha_2}-\alpha_2 \beta_2
v^{\alpha_2}+(\alpha_3-1) \alpha_3 \beta_3 v^{\alpha_3})
(\alpha_2 \beta_2 x v^{\alpha_2-3}-\alpha_4 \beta_4 y
v^{\alpha_4-3})$,
\item
$D'(q)=
\beta_2 v^{\alpha_2-2} (\alpha_2 \beta_2 x
v^{\alpha_2}+\alpha_3 \beta_3 x v^{\alpha_3}+\delta_3 v
(x+y)-\delta_4 v-v x^2-2 v x y-v y^2) (\alpha_2 \beta_3 x
v^{\alpha_3}-\alpha_2 \beta_4 x v^{\alpha_4}+(\alpha_2-1) \delta_3
v (x+y)+\delta_4 (v-\alpha_2 v)-\alpha_2 v x^2-2 \alpha_2 v x
y-\alpha_2 v y^2-\alpha_3 \beta_3 x v^{\alpha_3}+\alpha_4 \beta_4
x v^{\alpha_4}+v x^2+2 v x y+v y^2)
+(\alpha_2 \beta_2 x
v^{\alpha_2-1}-\alpha_4 \beta_4 y v^{\alpha_4-1})
(( \alpha_2-1) \beta_2 v^{\alpha_2} (\delta_4-(x+y)
(\delta_3-x-y))+(\alpha_3-1) \beta_3 v^{\alpha_3} (y
(y-\delta_3)+\delta_4-x^2)+(\alpha_4-1) \beta_4 y v^{\alpha_4}
(\delta_3-2 (x+y)))
$,
\item
$
D''(q)=
(\beta_2 x v^{\alpha_2-1}-\beta_4 y v^{\alpha_4-1})
(\alpha_2 \beta_2 x v^{\alpha_2}+\alpha_3 \beta_3 x
v^{\alpha_3}+\delta_3 v (x+y)-\delta_4 v-v x^2-2 v x y-v
y^2)-(\alpha_2 \beta_2 x v^{\alpha_2-1}-\alpha_4 \beta_4 y
v^{\alpha_4-1}) (\beta_2 x v^{\alpha_2}+\beta_3 x
v^{\alpha_3}+\delta_3 v (x+y)-\delta_4 v-v x^2-2 v x y-v y^2)
$,
\end{itemize}
and the singular control is given by: $u_s = -D'(q) / D(q)$.

\subsubsection{Classification of local syntheses for the McKeithan network}

We consider the case ${\max\ [A]}$, with $x=[A]$. We proceed as follows.

{\bf Stratification of the terminal manifold: $x=d$.}

\paragraph*{Singular locus} ${\cal S}: n\cdot [G,F](q)=0 \text{ and } x=d$ with $n=(1,0,0)$.
It is given by:
\begin{equation}
\begin{aligned}
  {\cal S}:\,  \alpha_2 \beta_2 d v^{\alpha_2-1}+\alpha_3 &\beta_3 d
      v^{\alpha_3-1}  +d \delta_3-\delta_4
      \\ 
      &-d^2+y (\delta_3-2 d)-y^2=0.
\end{aligned}
      \label{eq:eqS}
    \end{equation}
Denoting by $\Delta$ the discriminant of the polynomial function 
$y\mapsto n\cdot [G,F](q) \cap x=d$, a singularity can occur for $\Delta=0$.

One has
\begin{lemma}
  Assume $\alpha_i , \beta_i, \delta_i >0$,\,i=1,2 and $d,v>0$. Then we have
  $\Delta=(\delta_1-\delta_2)^2+4d\,( \alpha_2\beta_2v^{\alpha_2-1}+\alpha_3\beta_3v^{\alpha_3-1} )> 0$ so that there is no ramification and ${\cal S}$ contains at most two real positive branches.
\end{lemma}

\begin{definition}
  A semi-bridge occurs at a point $q\in {\cal S}$ if $n\cdot \left[ \left[ G,F \right],G \right](q)=0$.
\end{definition}

Computing, a semi-bridge occurs if 

\begin{equation}
\label{cond:semi-bridge}
v^{\alpha_3-\alpha_2}=-\frac{(\alpha_2-1)\alpha_2\beta_2}{(\alpha_3-1)\alpha_3\beta_3}.
\end{equation}

\paragraph*{Exceptional locus}
It is given by ${\cal E}: n\cdot F(q)=0$ and $x=d$.
Computing, one gets:
\begin{equation}
\begin{aligned}
{\cal E}:\,
   -\beta_2 d v^{\alpha_2}&-\beta_3 d v^{\alpha_3}+d^2 v-d \delta_3 v
   \\
   &+y\, (2 d v-\delta_3 v)+\delta_4 v+v y^2 = 0.
       \end{aligned}
  \label{eq:eqE}
\end{equation}
The discriminant of the polynomial $n\cdot X(q)$ in $y$ is
$\Delta =    v (4 d (\beta_2 v^{\alpha_2}+\beta_3$ $
    v^{\alpha_3})+v(\delta_1-\delta_2)^2)>0$
    and ${\cal E}$ contains at most two real positive branches.

\begin{figure}
    \centering
    \def\svgwidth{0.47\textwidth}
    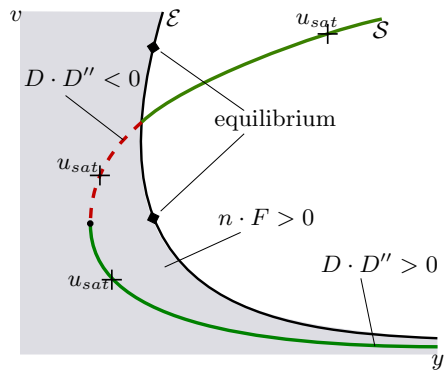
  \caption{Stratification of the surface $x=d$ for the McKeithan reaction.
  Dotted line: elliptic, red line: hyperbolic.}
  \label{fig:stratification_mck}
\end{figure}

Fig.\ref{fig:stratification_mck} gives a picture of stratification
of $x=d$ for the McKeithan system with a focus where 
${\cal S}$ is folded.
We represent the sets ${\cal S}$ and ${\cal E}$ and 
the stratification of ${\cal S}$ in 
hyperbolic, elliptic and parabolic points.
From this, we deduce the admissible points.
We are in the flat case and a point of $N$ is either
an ordinary point or a fold point.
For an ordinary point, the final optimal control is regular and is equal
to $-\textrm{sign}\ \Phi(0)$. At a fold point on ${\cal S}$, the final
optimal control may be singular (see Definition \ref{def:fold_case}).  

The optimal syntheses near the fold can be described using the techniques of the tutorial
model \eqref{eq:normal-form} and symbolic computations.

\section{Conclusion}
  \label{sec5}
  In this article, we have presented the general tools to analyze a Mayer problem
  to optimize the yield of chemical networks. Using the McKeithan network, we have 
  completed the analysis of \cite{bonnard1998} for the simple scheme
  $A\rightarrow B \rightarrow C$.
  
  The main point is to extend a geometric techniques from \cite{bonnard1997} 
  to consider the case
  when the strict Legendre-Clebsch is not satisfied. 
  The starting point is to analyze the semi-bridge phenomenon.
  
  Note that with the same techniques we can analyze the situations
  up to codimension $2$ describe in details in \cite{launay1997} 
  (in particular near the exceptional locus ${\cal E}$).
  Furthermore, it can be extended to an important problem
	of codimension $3$ occurring for reversible 
  chemical networks and the existence of equilibria, see Fig. \ref{fig:stratification_mck}.

\section{Appendix}

\subsection{Strata.m}
\label{prog:strata}
\begin{lstlisting}
(* Strata.m -- Tested with Mathematica v11.3 Linux *)
(*
This program computes the singular set, switching surface \
and splitting locus to derive the time minimal synthesis in a neighborhood of N.
*)
(* Lie derivative of w wrt v *) 
LieDerive[v_List,w_List,var_List]:=Transpose[D[w,#]& /@ var] . v /; \
	Length[v] == Length[var]

(* Lie bracket [v,w] *)
LieBracket[v_List,w_List,var_List]:=\
	LieDerive[w,v,var] - LieDerive[v,w,var] /; \
	Length[v] == Length[w] == Length[var]

(* Poisson bracket {f,g} *)
PoissonBracket[f_, g_,q_List,p_List] /; Length[q] == Length[p] := \
	Total @ Flatten @ Fold[List,0,MapThread[D[f,#1] D[g,#2] - D[f,#2] D[g,#1] &, {q, p}]]

(* Remove monomials of order >= ord *)
EliminTe[expr_,var_List,ord_Integer]:=Module[{},\
	FromCoefficientRules[Select[CoefficientRules[expr,var],Total@#[[1]] <= ord &], var]]

(* Multivariate Taylor expansions *)
mTaylor[expr_,var_List,pts_List,ord_Integer] /; \
	Length[var] == Length[pts] := Normal[Series[(expr /. Thread[var -> s (var - pts) + pts])//ExpandAll, {s, 0, ord}]] /. {s -> 1}

(* 
	Integration in (q,p) of H=p.(XX + u YY) for small time in o(t^ord) 
*)
IntSmallTime[XX_,YY_,u_,qfpf_List,ord_Integer]:=Module[{res},
  pv          = {p1,p2,p3};
  qv          = {x,y,z};
  H           = {p1,p2,p3}.(XX + u YY);
  res         = Sum[t^k / k! Nest[PoissonBracket[#,H,qv,pv]&,#,k],{k,0,ord}] \
                /. Thread[Join[qv,pv]->qfpf]& /@ Join[qv,pv]
];


$Assumptions = (ep==1 || ep==-1);

(*--------------------------------------------------------*)
(*--------------------------main--------------------------*)
(*--------------------------------------------------------*)

(*.......... constants ..........*)
a1      		= 1;      a2 = 1;      a4 = 1;
b       		= 1;       c = 1;       d = 0;	
uh0     		= 1;      uhx = 1;    uhy = 1;

(*.......... model ..........*)
XX      		= {1 + a z^2 + a2 x z^2 + a3 x y^2 + a4 y z^2,\
				b z, -uh0 -uhx x -uhy y + c z};
X       		= XX /. {x-> x-d};
Y       		= {0,0,1};

(*.......... Lie brackets computations ..........*)
YX      		= LieBracket[Y,X,{x,y,v}];
YXY     		= LieBracket[YX,Y,{x,y,v}];
YXX     		= LieBracket[YX,X,{x,y,v}];
DD      		= Det[{Y,YX,YXY}];
DDp     		= Det[{Y,YX,YXX}];
DDpp    		= Det[{Y,YX,X}];
aux     		= Solve[(YX[[1]])==0,{y}];
Sy[v_]  		= Expand[(y /. aux[[1]]) /. x->0]

(*.......... Gammas ..........*)
ord         = 3;
us          = - DDp / DD;
q0p0        = {0,w0,0,1,0,0};
qpt         = IntSmallTime[X,Y,us,q0p0,ord];
ps          = {t0val,w0val,s0val};
gms         = mTaylor[#,{t,w0,s0},ps,ordG]& /@ qpt[[1;;3]] //Expand
gmsall      = mTaylor[#,{t,w0,s0},ps,ordG]& /@ qpt //Expand

(*.......... Gamma+,Gamma- ..........*)
u           = ep;
q0p0        = {0,w0,s0,1,0,0};
pt          = {t0val,w0val,s0val};
qpt         = Refine /@ IntSmallTime[X,Y,u,q0p0,ord];
qpt         = mTaylor[#,{t,w0,s0},pt,ordG] & /@ qpt;
gmall       = Refine /@ qpt //Refine
gmp         = gmall[[1;;3]] /. {ep->1}
gmm         = gmall[[1;;3]] /. {ep->-1}


(*.......... W+,W- ..........*)
ord         = 3;
qpt         = Refine /@ IntSmallTime[X,Y,u,q0p0,ord];
eq          = mTaylor[qpt[[6]],vars,{t0val,w0val,s0val},ord] // Simplify
aux         = Solve[eq==0,vw];
sols        = Flatten[vw /. aux];
sols        = Refine[Expand[mTaylor[#,{t,w0,s0},{0,w0val,s0val},ord] & /@ sols]] // Simplify
solp        = mTaylor[ # /. {ep->1},vars,{t0val,w0val,s0val},ordG] & /@ sols;
solm        = mTaylor[ # /.{ep->-1},vars,{t0val,w0val,s0val},ordG] & /@ sols
qptp        = (qpt /. ep->1 /. vw-># ) & /@ solp;
qpwps       = mTaylor[#,{t,w0,s0},pt,ordG] & /@ qptp // Flatten;
qptm        = (qpt /. ep->-1 /. vw-># ) & /@ solm;
qpwms       = mTaylor[#,{t,w0,s0},pt,ordG] & /@ qptm //Flatten

(*.......... Ws ..........*)
ord         = 3;
q0p0        = {0,w0,0,1,0,0};
qpt         = Refine /@ IntSmallTime[X,Y,us,q0p0,ord];
qpt         = mTaylor[#,{t,w0,s0},pt,ord] & /@ qpt;
q0p02       = qpt /. t->ss;
qpt2        = Refine /@ IntSmallTime[X,Y,1,q0p02,ord];
qpt2        = mTaylor[#,{t,w0,s0,ss},{0,w0val,s0val,0},ord] & /@ qpt2;
p3t2        = qpt2[[6]]
aux         = Solve[p3t2==0,ss];
sols        = Flatten[ss /. aux];
sols        = Refine[Expand[mTaylor[#,{t,w0,s0},{0,w0val,s0val},ord] & /@ sols]] 
qpws        = qpt2 /. {ss-> #} & /@ sols;
qpws        = mTaylor[#,{t,w0,s0},pt,ordG] & /@ qpws //Flatten
ws          = Refine /@ qpws[[1;;3]]

(*.......... C1 ..........*)
ord         = 3;
q0p0        = {0,w0,s0,1,0,0};
pt          = {0,w0val,s0val};
qpt         = IntSmallTime[X,Y,u,q0p0,ord];
qpt         = mTaylor[#,{t,w0,s0},pt,ord] & /@ qpt;
(* u=-1 *)
{xp,yp,zp}  = qpt[[1;;3]] /.{ep->-1};
(* u=+1 *)
{xm,ym,zm}  = qpt[[1;;3]] /.{ep->1,s0->s0p,w0->w0p};

eqns        = {xp==xm,yp==ym,zp==zm}; 
aux         = Solve[eqns,{t,w0p,s0p}] // Expand;
aux2        = {t,w0p,s0p} /. aux;
c1sols      = mTaylor[#,{t,w0,s0},pt,ord] & /@ aux2;

(*.......... C12 ..........*)
ord			= 2;
q0p0			= {0,w0,s0,1,0,0};
pt				= {0,w0val,s0val};
qpt			= IntSmallTime[X,Y,u,q0p0,ord];
qpt			= mTaylor[#,{t,w0,s0},pt,ord] & /@ qpt;
vars			= {t,ss1,ss2,w0,s0,w0p,s0p};
pt2 			= {0,0,0,w0val,s0val,w0val,s0val};

(* u=-1 *)
{x1,y1,z1}  = qpt[[1;;3]] /.{ep->-1};

(* u=-1 & switch *)
ptfin       = (qpt /.{ep->-1,w0->w0p,s0->s0p})/.t->ss1;
qpt2        = IntSmallTime[X,Y,1,ptfin,ord] /. t->ss2;
{x3,y3,z3}  = mTaylor[#,vars,pt2,ord] & /@ qpt2[[1;;3]];

eqns        = {x1==x3,y1==y3,z1==z3,ss2+ss1==t,ptfin[[6]]==0};
aux         = Solve[eqns,{t,ss1,ss2,w0p,s0p}];
aux2        = {t,ss1,ss2,w0p,s0p} /. aux;
c12sols     = mTaylor[#,vars,pt2,ord] & /@ sols;


\end{lstlisting}

\end{document}